\documentclass[pdflatex,12pt,sn-mathphys-num]{sn-jnl}

\usepackage{amsmath,amssymb,amsfonts,bm}
\usepackage{graphicx}
\usepackage{multirow}
\usepackage{xcolor}
\usepackage{url}
\usepackage[stable]{footmisc}

\newtheorem{remark}{Remark}
\newtheorem{proposition}{Proposition}

\usepackage{algorithm}
\usepackage{algpseudocode}

\newcommand{\norm}[1]{\left\lVert #1\right\rVert}

\newcommand{\normA}[1]{\left\lVert #1\right\rVert_{\mathbf{A}}}
\newcommand{\lammin}{\lambda_{m}}
\newcommand{\lammax}{\lambda_{M}}

\begin{document}

\title{Scaled Coordinate Descent for Convex Quadratic Programs}

\author*[1]{\fnm{Cristian} \sur{Rusu}}
\email{cristian.rusu@unibuc.ro}
\affil[1]{\orgdiv{Faculty of Mathematics and Computer Science},
            \orgname{University of Bucharest},
            \orgaddress{\city{Bucharest},
            \postcode{030018},
            \country{Romania}}}

\abstract{
We study a modification of coordinate descent for unconstrained symmetric positive definite quadratic programs that combines a coordinate correction with an exact scaling of the current iterate. For a selected coordinate, the method minimizes the objective function over the two-dimensional subspace spanned by the current iterate and the corresponding standard basis vector. Closed-form expressions are obtained for the two update parameters and for the resulting objective decrease. The scaling produces a state-dependent effective curvature that is no larger than the diagonal curvature used by standard coordinate descent, so the exact scaled update makes at least as much progress for the same selected coordinate. We discuss initialization, numerical stability, randomized coordinate selection, experimental comparisons with ordinary two-coordinate descent, and a block coordinate extension. In the spirit of previous theoretical results, we give a probabilistic linear convergence result. Finally, we provide experimental results to show the effectiveness of the proposed method on synthetic data.
}


\maketitle

\section{Introduction}

Coordinate Descent (CD) methods~\cite{Wright2015} form a broad class of optimization techniques that are particularly useful in signal processing~\cite{9420314,10736441} and machine learning~\cite{JSSv033i01,pmlr-v97-zeng19a} applications. Their popularity follows from their conceptual simplicity, low cost per iteration, suitability for large-scale problems, and availability in software packages~\cite{glmnet,10.5555/1953048.2078195}. At each iteration, a CD method updates one variable, or a small block of variables, while keeping the remaining fixed. Classical convergence results for convex differentiable objectives were established in early work such as~\cite{LuoTseng1992}, while modern treatments cover cyclic, randomized, greedy, accelerated, parallel, block, and proximal variants~\cite{Wright2015}.

Given a symmetric positive definite matrix $\mathbf{A}\in\mathbb{R}^{n\times n}$ and a vector $\mathbf{b}\in\mathbb{R}^n$, we consider the strictly convex quadratic problem
\begin{equation}
    \underset{\mathbf{x}\in\mathbb{R}^n}{\operatorname{minimize}}
    \
    f(\mathbf{x})
    =\frac{1}{2}\mathbf{x}^{T}\mathbf{A}\mathbf{x}
    -\mathbf{b}^{T}\mathbf{x}.
    \label{eq:intro-quadratic-problem}
\end{equation}
The unique minimizer of~\eqref{eq:intro-quadratic-problem} is the solution of $\mathbf{A}\mathbf{x}=\mathbf{b}$. A standard exact coordinate step selects an index $i$ and solves the one-dimensional problem
\begin{equation}
    \underset{\alpha\in\mathbb{R}}{\operatorname{minimize}}
    \ f\bigl(\mathbf{x} + \alpha\mathbf{e}_i \bigr),
    \label{eq:intro-standard-coordinate-step}
\end{equation}
where $\mathbf{e}_i$ is the $i$th canonical basis vector of $\mathbb{R}^n$. For a quadratic objective, this update is algebraically identical to a Gauss--Seidel relaxation step for the linear system $\mathbf{A}\mathbf{x}=\mathbf{b}$. Coordinate descent for a full column-rank least squares problem can similarly be interpreted as Gauss--Seidel applied to the normal equations, although practical implementations generally update the least squares residual without explicitly forming the normal matrix. This connection has motivated extensive work on randomized Gauss--Seidel and coordinate methods for linear systems and least squares problems~\cite{LeventhalLewis2010,Dumitrescu2015,doi:10.1137/15M1014425}. Randomized Kaczmarz provides the complementary row-action viewpoint~\cite{StrohmerVershynin2009}, while the sketch-and-project framework places coordinate descent, Kaczmarz, and several block variants within a common projection formalism~\cite{GowerRichtarik2015}.

Randomized coordinate selection is attractive because it yields global expected convergence estimates in terms of spectral or coordinate-wise condition measures. Sampling coordinates in proportion to their Lipschitz constants is standard for smooth quadratic objectives and underlies the complexity analyses of randomized and accelerated coordinate methods~\cite{doi:10.1137/100802001,Lu2015,RichtarikTakac2014}. Accelerated parallel and proximal variants improve worst-case iteration bounds while retaining inexpensive partial updates~\cite{doi:10.1137/130949993}. A complementary line of work seeks greater progress per iteration through greedy selection. Gauss--Southwell-type rules choose coordinates with large partial derivatives~\cite{Nutini2015}, while greedy randomized and relaxed greedy strategies use residual information to construct an active set of promising coordinates~\cite{BaiWu2019,ZhangGuo2020}. More general analyses separate the effects of the block partition, the block-selection rule, and the local block solver~\cite{Nutini2022}. These approaches alter how the coordinate or block is selected, but the local correction is still performed in a fixed coordinate subspace.

Block coordinate descent enlarges the search line in~\eqref{eq:intro-standard-coordinate-step} to a prescribed coordinate subspace and solves the corresponding local problem either exactly or approximately. Classical convergence results include~\cite{Tseng2001,BeckTetruashvili2013,RichtarikTakac2014}, while subspace-correction theory provides a broader numerical interpretation based on a fixed decomposition of the ambient space~\cite{Xu1992}. Sequential subspace optimization further accelerates coordinate and first-order methods by minimizing over a small space containing several algorithmically generated directions~\cite{EladMatalonZibulevsky2007}. For quadratic objectives, the progress of an exact coordinate-block correction depends on the inverse of the associated principal submatrix. This observation underlies determinant-based volume sampling, which selects a coordinate set $S$ with probability proportional to $\det(\mathbf{A}_{SS})$ and obtains convergence rates that improve with the sampled block size~\cite{RodomanovKropotov2020}.

The methods most closely related to the comparisons developed in this paper use two data-dependent directions per iteration. The two-subspace Kaczmarz method combines two row projections and is particularly effective for coherent systems~\cite{NeedellWard2013}. For column-action least squares solvers, two-step and multi-column Gauss--Seidel variants select and jointly process multiple columns~\cite{LiuJiangGu2021,LiaoLuYin2022}. Of particular relevance, greedy double-subspace coordinate descent selects two active columns and uses Gram--Schmidt orthogonalization to form a two-dimensional correction, with its strongest advantages arising for highly coherent matrices~\cite{JinLi2023}. Related two-subspace and extended methods include~\cite{Wu2022,MustafaSaha2024,MustafaSaha2026,ZhangDaiZhao2026}. Although these methods differ in their selection rules and local solvers, they all enlarge the local space using two rows, two columns, or a fixed coordinate block.

Low-dimensional acceleration methods provide another relevant perspective. Anderson acceleration forms an extrapolated iterate from a short history of fixed-point residuals and has also been adapted to coordinate descent~\cite{10.1145/3197517.3201290,BertrandMassias2021}. Coordinate descent has been combined with Conjugate Gradient (CG) acceleration for linear systems~\cite{Gordon2023}, and two-dimensional minimum-residual corrections have been proposed for accelerating stationary iterations~\cite{BeikBenziNajafi2024}. Such techniques enrich the search space with previous steps, residuals, or conjugate directions. By contrast, the method considered here stores no history window and constructs no accumulated Krylov basis.

In this paper, we study a simple modification of the exact coordinate update. In iteration $k$, selecting the coordinate $i$, the next iterate is obtained:
\begin{equation}
    \mathbf{x}^{(k+1)}
    =\beta^\star\mathbf{x}^{(k)}+\alpha^\star\mathbf{e}_i,
    \text{ and }
    (\beta^\star,\alpha^\star) = 
    \underset{\beta,\alpha\in\mathbb{R}}{\operatorname{arg\,min}}
    f\bigl(\beta\mathbf{x}^{(k)}+\alpha\mathbf{e}_i\bigr).
    \label{eq:intro-scaled-coordinate-step}
\end{equation}
Thus, the current iterate is optimally rescaled while the selected coordinate is corrected, and the objective is minimized over the state-dependent plane
$\operatorname{span}\{\mathbf{x}^{(k)},\mathbf{e}_i\}$. The noncoordinate direction is not selected from a fixed dictionary: it is the current, generally dense, iterate and therefore carries information accumulated over all preceding updates. The two parameters are chosen jointly by minimizing the objective over the span of the current iterate and the selected coordinate vector. We derive the exact update, study its effective curvature and numerical stability, give an adaptive randomized sampling rule and its expected convergence estimate, compare the method with ordinary two-coordinate descent and volume sampling, and extend the construction to block-coordinate updates.

\section{Scaled coordinate descent for convex quadratic programs}

For the problem in \eqref{eq:intro-standard-coordinate-step} we know that the optimal solution is given by:
\begin{equation}
	\alpha^\star  =  - \frac{\mathbf{d}^T(\mathbf{Ax} - \mathbf{b})}{\mathbf{d}^T \mathbf{Ad}},
	\label{eq:best_alpha_original}
\end{equation}
and the objective function value is:
\begin{equation}
    \min_{\alpha} f(\mathbf{x}+\alpha\mathbf{d})
    =f(\mathbf{x})-
    \frac{1}{2}
    \frac{\bigl(\mathbf{d}^T(\mathbf{Ax}-\mathbf{b})\bigr)^2}
    {\mathbf{d}^T\mathbf{Ad}}.
    \label{eq:1d_search}
\end{equation}

In this paper, we analyze a slight modification of the line search optimization problem formulated above. At each step of the algorithm, we consider a given search direction $\mathbf{d} \in \mathbb{R}^n, \mathbf{d} \neq \mathbf{0}_{n \times 1},$ and we want to perform exact line and scaling updates such that we solve:
\begin{equation}
	\underset{\alpha, \beta \in \mathbb{R}}{\text{minimize}} \ f(\beta \mathbf{x} + \alpha \mathbf{d}).
	\label{eq:qp_2d}
\end{equation}

The line search is performed along $\mathbf{d}$ and the scaling update is applied to $\mathbf{x}$. This can be seen as an overall scaling and translation along the search direction. For $\beta = 1$ we recover the ordinary one-dimensional coordinate descent problem previously described. In the following subsection, we explore the computational solution and theoretical insights to the problem in \eqref{eq:qp_2d}.

\subsection{Preliminaries on the scaling factor optimization}

First, we isolate the effect of scaling the current vector in the quadratic optimization problem. Then we have the following remark.

\begin{remark}[Optimal scaling]
For a fixed nonzero vector $\mathbf{x}$, consider:
\begin{equation}
    \underset{\beta\in\mathbb{R}}{\text{minimize}}\
    f(\beta\mathbf{x})
    =\frac{\beta^2}{2}\mathbf{x}^T\mathbf{Ax}
    -\beta\mathbf{b}^T\mathbf{x}.
    \label{eq:optimally_scaled}
\end{equation}
Then the following holds:
\begin{enumerate}
    \item The optimal scaling is: $\beta^\star=\frac{\mathbf{b}^T\mathbf{x}}
        {\mathbf{x}^T\mathbf{Ax}}$;

    \item The minimum objective value is: $f(\beta^\star\mathbf{x})
        = -\frac{1}{2}\frac{(\mathbf{b}^T\mathbf{x})^2}
        {\mathbf{x}^T\mathbf{Ax}}$;

    \item If $\mathbf{x} = \mathbf{b}$, then $\beta^\star=\frac{\|\mathbf{b}\|_2^2}{\mathbf{b}^T\mathbf{Ab}}$, $f(\beta^\star\mathbf{b})
        =-\frac{1}{2}
        \frac{\|\mathbf{b}\|_2^4}
        {\mathbf{b}^T\mathbf{Ab}}$,
    and we have: $-\frac{\|\mathbf{b}\|_2^2}{2\lambda_m(\mathbf{A})}
        \leq f(\beta^\star\mathbf{b})
        \leq -\frac{\|\mathbf{b}\|_2^2}{2\lambda_M(\mathbf{A})}$. At $\mathbf{x}^\star=\mathbf{A}^{-1}\mathbf{b}$ we have $f(\mathbf{x}^\star) =-\frac{1}{2}\mathbf{b}^T\mathbf{A}^{-1}\mathbf{b}$.
\end{enumerate}
\end{remark}
\noindent\textbf{Proof.}
These statements follow directly by minimizing the scalar quadratic in $\beta$ and the quadratic form eigenvalue bounds. \hfill$\blacksquare$

For any $\mathbf{x}$ that is optimally scaled in \eqref{eq:optimally_scaled} the equality:
\begin{equation}
    \mathbf{b}^T \mathbf{x} = \mathbf{x}^T \mathbf{Ax},
    \label{eq:scaling_invariant}
\end{equation}
holds and is maintained. Equivalently, with the residual vector:
\begin{equation}
\mathbf{r} = \mathbf{Ax} - \mathbf{b},
\end{equation}
this gives the orthogonality property $\mathbf{x}^T\mathbf{r} = 0$.
%The choice of a vector proportional to $\mathbf{b}$ is therefore natural for initialization.%The zero vector is less useful because it does not provide a nonzero scaling direction.

\subsection{The solution to the proposed scaled descent method}

Given the effect of the scaling variable $\beta$, we propose to jointly optimize it and the search direction, leading to the following closed-form result.

\begin{proposition}[Optimal joint scaling and coordinate descent]
Assume that $\mathbf{x}$ and $\mathbf{d}$ are linearly independent. The solution of~\eqref{eq:qp_2d} is:
\begin{equation}
    \beta^\star =
    \frac{(\mathbf{d}^T\mathbf{Ad})(\mathbf{b}^T\mathbf{x})
    -(\mathbf{d}^T\mathbf{Ax})(\mathbf{b}^T\mathbf{d})}
    {(\mathbf{d}^T\mathbf{Ad})(\mathbf{x}^T\mathbf{Ax})
    -(\mathbf{d}^T\mathbf{Ax})^2},
    \alpha^\star =-
    \frac{(\mathbf{d}^T\mathbf{Ax})(\mathbf{b}^T\mathbf{x})
    -(\mathbf{x}^T\mathbf{Ax})(\mathbf{b}^T\mathbf{d})}
    {(\mathbf{d}^T\mathbf{Ad})(\mathbf{x}^T\mathbf{Ax})
    -(\mathbf{d}^T\mathbf{Ax})^2}.
\label{eq:sdp_optimal_alpha_beta}
\end{equation}
\end{proposition}
\noindent\textbf{Proof.}
The problem is the two-variable quadratic:
\begin{equation}
    \underset{\alpha,\beta\in\mathbb{R}}{\text{minimize}}\
    \frac{1}{2}
    \begin{bmatrix}\alpha&\beta\end{bmatrix}
    \begin{bmatrix}
        \mathbf{d}^T\mathbf{Ad} & \mathbf{d}^T\mathbf{Ax}\\
        \mathbf{d}^T\mathbf{Ax} & \mathbf{x}^T\mathbf{Ax}
    \end{bmatrix}
    \begin{bmatrix}\alpha\\\beta\end{bmatrix}-
    \begin{bmatrix}\mathbf{b}^T\mathbf{d}&\mathbf{b}^T\mathbf{x}\end{bmatrix}
    \begin{bmatrix}\alpha\\\beta\end{bmatrix}.
\label{eq:the_detailed_problem}
\end{equation}
The matrix in~\eqref{eq:the_detailed_problem} is positive definite when the two directions are linearly independent, and solving the corresponding $2\times2$ linear system gives~\eqref{eq:sdp_optimal_alpha_beta}. \hfill$\blacksquare$

Assuming~\eqref{eq:scaling_invariant} holds before the update, then:
\begin{equation}
    \beta^\star = \frac{(\mathbf{d}^T \mathbf{Ad})(\mathbf{b}^T \mathbf{x})
    - (\mathbf{d}^T \mathbf{Ax})(\mathbf{b}^T \mathbf{d})}{(\mathbf{d}^T\mathbf{Ad})(\mathbf{b}^T\mathbf{x})
    -(\mathbf{d}^T\mathbf{Ax})^2},
    \alpha^\star
     =- \frac{\mathbf{d}^T(\mathbf{Ax}-\mathbf{b}) (\mathbf{b}^T\mathbf{x})} {(\mathbf{d}^T\mathbf{Ad})(\mathbf{b}^T\mathbf{x})
    -(\mathbf{d}^T\mathbf{Ax})^2},
\label{eq:sdp_optimal_alpha_beta_simplified}
\end{equation}
and the objective value is
\begin{equation}
    f(\beta^\star\mathbf{x}+\alpha^\star\mathbf{d})
    =f(\mathbf{x})-
    \frac{1}{2}
    \frac{(\mathbf{b}^T\mathbf{x})
    \bigl(\mathbf{d}^T(\mathbf{Ax}-\mathbf{b})\bigr)^2}
    {(\mathbf{d}^T\mathbf{Ad})(\mathbf{b}^T\mathbf{x})
    -(\mathbf{d}^T\mathbf{Ax})^2}.
    \label{eq:general_scaled_decrease}
\end{equation}

The adaptation to a coordinate descent method follows immediately.

\begin{proposition}[Optimal coordinate scaling formulas]
Let the search direction be fixed to $\mathbf{d} = \mathbf{e}_i$, and let $\mathbf{a}_i$ denote the $i$th column of $\mathbf{A}$. Under~\eqref{eq:scaling_invariant}, the optimal parameters are:
\begin{equation}
    \beta^\star
    =\frac{A_{ii}(\mathbf{x}^T\mathbf{Ax})
    -b_i(\mathbf{a}_i^T\mathbf{x})}
    {A_{ii}(\mathbf{x}^T\mathbf{Ax})
    -(\mathbf{a}_i^T\mathbf{x})^2} = 1+
    \frac{(\mathbf{a}_i^T\mathbf{x})
    (\mathbf{a}_i^T\mathbf{x}-b_i)}
    {A_{ii}(\mathbf{x}^T\mathbf{Ax})
    -(\mathbf{a}_i^T\mathbf{x})^2},
\label{eq:sdp_best_beta}
\end{equation}
\begin{equation}
    \alpha^\star
    =-
    \frac{(\mathbf{x}^T\mathbf{Ax})
    (\mathbf{a}_i^T\mathbf{x}-b_i)}
    {A_{ii}(\mathbf{x}^T\mathbf{Ax})
    -(\mathbf{a}_i^T\mathbf{x})^2}.
    \label{eq:sdp_best_alpha}
\end{equation}
The corresponding objective value is:
\begin{equation}
    f(\beta^\star \mathbf{x} + \alpha^\star \mathbf{e}_i) = f(\mathbf{x}) -
    \frac{1}{2}
     \frac{(\mathbf{a}_i^T \mathbf{x}-b_i)^2}
    {A_{ii}-\dfrac{(\mathbf{a}_i^T \mathbf{x})^2}
    {\mathbf{x}^T \mathbf{Ax}}} = f(\mathbf{x})-
    \frac{1}{2}
    \frac{(\mathbf{a}_i^T \mathbf{x}-b_i)^2}
    {A_{ii}(1-\gamma_i)},
\label{eq:new_objective}
\end{equation}
where
\begin{equation}
    \gamma_i=
    \frac{(\mathbf{e}_i^T\mathbf{Ax})^2}
    {A_{ii}(\mathbf{x}^T\mathbf{Ax})}.
    \label{eq:gammai}
\end{equation}
The Cauchy--Schwarz inequality in the $\mathbf{A}$-inner product gives $0\leq\gamma_i\leq1$. If $\mathbf{x}$ is not proportional to $\mathbf{e}_i$, then $\gamma_i<1$. If $\mathbf{x}=c\mathbf{e}_i$ and~\eqref{eq:scaling_invariant} holds, then the denominator and the selected residual both vanish. \hfill$\blacksquare$
\end{proposition}

To highlight the difference between the proposed scaled method and the ordinary coordinate descent method, consider the following remark.

\begin{remark}[Comparison with ordinary coordinate descent]
For ordinary coordinate descent, the optimal update is:
\begin{equation}
\alpha^\star_\text{CD}
    = -\frac{\mathbf{a}_i^T \mathbf{x}-b_i}{A_{ii}}, \text{ and } \min_{\alpha} f(\mathbf{x} + \alpha\mathbf{e}_i)
    = f(\mathbf{x}) -
    \frac{1}{2} \frac{(\mathbf{a}_i^T \mathbf{x}-b_i)^2}{A_{ii}}.
\end{equation}
Since $A_{ii}(1-\gamma_i)\leq A_{ii}$, the exact scaled update makes at least as much objective decrease as ordinary coordinate descent for the same selected coordinate.
\end{remark}

Based on these updates, we propose the following scaled coordinate descent algorithm.

\subsection{The proposed algorithm}

The proposed Scaled Coordinate Descent (SCD) algorithm for convex quadratic programs is described in Algorithm~\ref{alg:scaled_coordinate_descent}. We discuss the main steps of the proposed algorithm and drop the index $k$, with the understanding that these insights hold at all steps.
\begin{algorithm}[!t]
\caption{Scaled Coordinate Descent (SCD)}
\label{alg:scaled_coordinate_descent}
\small
\begin{algorithmic}[1]
\Require $\mathbf{A} \in \mathbb{R}^{n\times n}$ symmetric positive definite, $\mathbf{b} \in \mathbb{R}^n$, $\mathbf{b} \neq \mathbf{0}$.
\Require $\mathbf{x}^{(0)} \neq \mathbf{0}$ with $\mathbf{b}^T \mathbf{x}^{(0)} \neq 0$, maximum iterations $K$, threshold $\kappa_{\max} \geq 1$.
\Ensure Approximate solution $\mathbf{x}^{(K)}$.

\State Compute quantity $\omega^{(0)}\gets\mathbf{b}^T\mathbf{x}^{(0)}$;
\State Compute optimal scaling $\displaystyle \beta^{(0)} \gets \frac{\omega^{(0)}}{(\mathbf{x}^{(0)})^T \mathbf{A} \mathbf{x}^{(0)}}$;
\State Apply optimal scaling $\mathbf{x}^{(0)} \gets \beta^{(0)} \mathbf{x}^{(0)}$ and update $\omega^{(0)} \gets \beta^{(0)} \omega^{(0)}$.

\For{$k=1,\ldots,K$}
    \State Choose $i\in\{1,\ldots,n\}$;
    \State Compute scalar product $\eta_i \gets \mathbf{a}_i^T \mathbf{x}^{(k-1)}$;
    %\If{$\displaystyle 1-\frac{\eta_i^2}{A_{ii}\omega^{(k-1)}} \geq \frac{4\kappa_{\max}}{(\kappa_{\max}+1)^2}$}
    \If{$\displaystyle \eta_i^2 \leq \left(\frac{\kappa_{\max} -1}{\kappa_{\max}+1}\right)^2 A_{ii}\omega^{(k-1)}$}
        \State Compute $\displaystyle \beta^{(k)} \gets \frac{A_{ii}\omega^{(k-1)} - b_i\eta_i}{A_{ii}\omega^{(k-1)} - \eta_i^2}$ and $\displaystyle \alpha^{(k)} \gets -\frac{\omega^{(k-1)}(\eta_i - b_i)}{A_{ii}\omega^{(k-1)} - \eta_i^2}$;
        \State Update $\mathbf{x}^{(k)} \gets \beta^{(k)} \mathbf{x}^{(k-1)} + \alpha^{(k)} \mathbf{e}_i$ and $\omega^{(k)} \gets \beta^{(k)} \omega^{(k-1)} + \alpha^{(k)} b_i$;
    \Else
    \State Compute the ordinary coordinate descent solution
    $\alpha_\text{CD}^{(k)}
    \gets
    -\frac{\eta_i-b_i}{A_{ii}}$;

    \State Compute the optimal scaling
    $\displaystyle
    \beta^{(k)}
    \gets
    \frac{
    \omega^{(k-1)}+\alpha_\text{CD}^{(k)}b_i
    }{
    \omega^{(k-1)}
    +\alpha_\text{CD}^{(k)}(\eta_i+b_i)
    }$;

    \State Update $\mathbf{x}^{(k)}
        \! \gets \!
        \beta^{(k)}
        (
        \mathbf{x}^{(k-1)}
        +\alpha_\text{CD}^{(k)}\mathbf{e}_i
        )$ and $\omega^{(k)}
        \gets
        \beta^{(k)}
        (
        \omega^{(k-1)}
        +\alpha_\text{CD}^{(k)}b_i
        )$;
    \EndIf
\EndFor
\State \textbf{return} $\mathbf{x}^{(K)}$.
\end{algorithmic}
\end{algorithm}

In the initialization step, given some random non-zero $\mathbf{x}^{(0)}$, we compute one matrix-vector product, costing $O(n^2)$ operations. This is the only step in the algorithm that has quadratic complexity. We could avoid this computational burden by initializing with the zero vector, i.e., $\mathbf{x}^{(0)} = \mathbf{0}_{n \times 1}$. The following remark is in order.

\begin{remark}[Initialization with the zero vector]
    If $\mathbf{x}^{(0)} = \mathbf{0}_{n \times 1}$, then a different initialization procedure is needed. At least the first iteration needs to perform an ordinary coordinate descent step and have $b_i \neq 0$, as the scaled direction does not contribute anything to the search space. Furthermore, supposing that the coordinates are visited in order during the first sweep, we see that before a coordinate $i$ is visited for the first time, the current iterate is Euclidean-orthogonal to $\mathbf{e}_i$. The Wielandt inequality then is $(\mathbf{x}^T\mathbf{Ae}_i)^2
    \leq
    \left(
    \frac{\kappa(\mathbf{A})-1}
    {\kappa(\mathbf{A})+1}
    \right)^2 \! \!
    (\mathbf{x}^T\mathbf{Ax})A_{ii}$,
and:
\begin{equation}
    A_{ii}-
    \frac{(\mathbf{x}^T\mathbf{Ae}_i)^2}
    {\mathbf{x}^T\mathbf{Ax}}
    \geq
    A_{ii}
    \left[
    1-
    \left(
    \frac{\kappa(\mathbf{A})-1}
    {\kappa(\mathbf{A})+1}
    \right)^2
    \right]
    =A_{ii}
    \frac{4\kappa(\mathbf{A})}
    {(\kappa(\mathbf{A})+1)^2} \approx \frac{4 A_{ii}}{\kappa(\mathbf{A})}.
\end{equation}
Thus, during the initial ordinary coordinate step, the denominators encountered during the remainder of the first sweep are guaranteed to be bounded away from zero. This is also related to the check on line 7 of Algorithm 1.
\end{remark}

%But, the purpose of the initialization process is to establish an initial estimate that is scaled optimally so that many of the desired properties and simplifications occur.

Then, in the iterative process, the scalar $\omega$ stores:
\begin{equation}
    \omega = \mathbf{b}^T \mathbf{x} = \mathbf{x}^T \mathbf{Ax} > 0.
\end{equation}
Consequently, the objective value at every scaled iterate is $f(\mathbf{x})=-\omega/2$. The exact update preserves this equality because it is the minimizer over $\operatorname{span}\{\mathbf{x},\mathbf{e}_i\}$. Note that in Algorithm~\ref{alg:scaled_coordinate_descent} line 7 we safeguard against a badly conditioned two-dimensional subspace basis. In the safeguarded branch, the algorithm first performs the ordinary coordinate descent correction
$\alpha_{\text{CD}}^\star = -(\eta_i-b_i)/A_{ii}$ and defines the corrected vector $\mathbf{x} + \alpha_{\text{CD}}^\star \mathbf{e}_i$. It then optimally scales the complete corrected vector, and its corresponding optimal scaling coefficient is:
\begin{equation}
    \beta^\dagger
    =
    \frac{
    \omega+\alpha_{\text{CD}}^\star b_i
    }{
    \omega
    +2\alpha_{\text{CD}}^\star\eta_i
    +(\alpha_{\text{CD}}^\star)^2A_{ii}
    }
    =
    \frac{
    \omega+\alpha_{\text{CD}}^\star b_i
    }{
    \omega+\alpha_{\text{CD}}^\star(\eta_i+b_i)
    }.
\end{equation}
The safeguarded iterate is therefore $\beta^\dagger
    \left(
    \mathbf{x}
    +\alpha_{\text{CD}}^\star \mathbf{e}_i
    \right)$,
and the maintained scalar is updated by $\beta^\dagger
    \left(
    \omega+\alpha_{\text{CD}}^\star b_i
    \right)$. By this construction, we maintain
$\mathbf{b}^{T}\mathbf{x} = \mathbf{x}^{T} \mathbf{A} \mathbf{x}
= \omega$, and hence
$\mathbf{x}^{T}
(\mathbf{A}\mathbf{x} - \mathbf{b}) = \mathbf{x}^{T} \mathbf{r} = 0$. Therefore, the safeguard preserves the optimal scaling invariant. The ordinary coordinate correction annihilates the selected residual
component before rescaling. After the final separate scaling, however, the $i^\text{th}$ component of the new residual is $r_i = (\beta^\dagger-1)b_i$ in general. Thus, the safeguarded branch preserves the scaling orthogonality but not necessarily the full two-dimensional Galerkin conditions of the exact scaled update. The exact convergence and expectation results stated below therefore refer specifically to the exact scaled update. Concerning the safeguard update, the following holds:
\begin{equation}
    f(\beta^\star \mathbf{x} + \alpha^\star \mathbf{e}_i) \leq f(\beta^\dagger
    \left(
    \mathbf{x}
    +\alpha_{\text{CD}}^\star\mathbf{e}_i
    \right)) \leq f(
    \mathbf{x}
    +\alpha_{\text{CD}}^\star\mathbf{e}_i) \leq f(\mathbf{x}).
\end{equation}

The safeguard conditioning check is based on the normalized Gram matrix of $\mathbf{x}$ and $\mathbf{e}_i$ in the $\mathbf{A}$-inner product. This is calculated from the normalization of the Gram matrix $\mathbf{G} = \begin{bmatrix} A_{ii} & \eta_i \\ \eta_i & \omega \end{bmatrix}$ by diagonal scaling to the equivalent matrix $\mathbf{\hat{G}} = \begin{bmatrix}
    1 & \sqrt{\gamma_i} \\
    \sqrt{\gamma_i} & 1
\end{bmatrix}$, where $\gamma_i$ is taken from \eqref{eq:gammai}. Its condition number is given by $\kappa(\mathbf{\hat{G}}) = \left( 1 + \dfrac{|\eta_i|}{\sqrt{A_{ii} \omega}} \right)
    \left( 1-\dfrac{|\eta_i|}{\sqrt{A_{ii} \omega}} \right)^{-1}$. Bounding this quantity to be smaller than $\kappa_{\max}$ leads to the condition used in Algorithm 1 at line 7.
% Unlike the condition number of the unscaled $2\times2$ Gram matrix, this quantity measures the near-linear dependence of the two search directions without being affected by a mismatch between their norms.

For a dense matrix, computing $\eta_i$ costs $O(n)$ operations, and scaling the iterate costs another $O(n)$ operations, for a total of $3n$ operations. Thus, the iteration cost remains linear in $n$ and is approximately one vector scaling more expensive than an ordinary coordinate update, which is $2n$ operations. The extra $n$ operations come from our new extra, but trivial, step of scaling the solution $\mathbf{x}$ by $\beta^\star$.

Note that we do not expand upon the choice of the coordinate descent direction in Line 5 of the proposed algorithm. Although index choice is very important, as we will also see in the following sections, Algorithm~\ref{alg:scaled_coordinate_descent} is generic enough to accommodate any selection rule: random, cyclic, greedy, etc.

The exact expectation formulas in Section~\ref{sec:expectation} apply to the exact update; the safeguard branch is included as a finite-precision implementation device.

\section{Theoretical aspects of the proposed scaling method}

In this section, we explore several theoretical aspects of the proposed scaling method. We are mainly focused on proving monotonicity, optimality and convergence rates. We also compare the proposed approach to the most natural reference, two-dimensional coordinate descent methods. While Algorithm~\ref{alg:scaled_coordinate_descent} uses a safeguarded update for badly conditioned practical scenarios, the results in this section refer only to the newly scaled updates developed in this paper.

%\subsection{Initialization}

%Algorithm~1 initializes the process with an optimally scaled nonzero vector, for example $\mathbf{b}$. Starting from the zero vector is also possible, but the first step must then be an ordinary one-dimensional coordinate update because the two-dimensional matrix in~\eqref{eq:the_detailed_problem} is singular when $\mathbf{x}=\mathbf{0}$.

\subsection{Monotonicity and local optimality}

As in ordinary coordinate descent, every exact scaled update given by \eqref{eq:new_objective} at step $k$ after choosing the coordinate descent direction $i$ is monotone, i.e., $f(\mathbf{x}^{(k+1)})\leq f(\mathbf{x}^{(k)}),$ because:
\begin{equation}
    f(\mathbf{x}^{(k)}) - f(\mathbf{x}^{(k+1)}) = \frac{1}{2} \frac{(r_i^{(k})^2}{ A_{ii} - \frac{ ( \mathbf{a}_i^T \mathbf{x}^{(k)} )^2 }{ (\mathbf{x}^{(k)})^T \mathbf{Ax}^{(k)}} }  = \frac{1}{2} \| \mathbf{x}^{(k+1)} - \mathbf{x}^{(k)} \|_\mathbf{A}^2.
    \label{eq:coordinate_decrease}
\end{equation}

For a nondegenerate update, $\alpha^\star=0$ exactly when the selected residual $r_i = \mathbf{a}_i^T\mathbf{x}-b_i$ is zero. Moreover, $\beta^\star=1$ when the selected residual is zero or when $\mathbf{a}_i^T\mathbf{x}=0$. In the latter case, the scaled update reduces to the ordinary one-dimensional coordinate update. From \eqref{eq:new_objective} we have that:
\begin{equation}
    0 < A_{ii}-
    \frac{(\mathbf{a}_i^T\mathbf{x})^2}
    {\mathbf{x}^T\mathbf{Ax}}
    \leq A_{ii},
    \label{eq:the_inequality}
\end{equation}
whenever the two search directions are independent. Hence, the exact scaled update makes at least as much progress as ordinary coordinate descent for the same coordinate.

We immediately have the following global convergence result.

\begin{proposition}[Global convergence]\label{thm:cyclic_convergence}
Assume there exists an integer $T\ge n$ such that every coordinate appears at least once in every $T$ consecutive iterations of the proposed scaling updates. Let the sequence $\{\mathbf{x}^{(k)}\}$ be generated by the exact scaled coordinate update, then:
\begin{equation}\label{eq:cyclic_convergence}
    \mathbf{x}^{(k)} \to \mathbf{x}^\star,
    \mathbf{r}^{(k)} \to \mathbf{0}_{n \times 1},
    f(\mathbf{x}^{(k)}) \to f(\mathbf{x}^\star).
\end{equation}
\end{proposition}
\noindent \textbf{Proof.} The objective is monotonically decreasing by \eqref{eq:coordinate_decrease} and bounded below by $f(\mathbf{x}^\star)$. Hence $\sum_{k=0}^{\infty} \bigl( f(\mathbf{x}^{(k)}) - f(\mathbf{x}^{(k+1)}) \bigr) < \infty$. By \eqref{eq:new_objective},
$\normA{\mathbf{x}^{(k+1)} - \mathbf{x}^{(k)}} \to 0$. Therefore, we have that:
\begin{equation}
\norm{\mathbf{r}^{(k+1)} - \mathbf{r}^{(k)}}_2
    = \norm{\mathbf{A}(\mathbf{x}^{(k+1)}-\mathbf{x}^{(k)})}_2
    \le \sqrt{\lammax(\mathbf{A})}\,
       \normA{\mathbf{x}^{(k+1)} - \mathbf{x}^{(k)}}
      \to 0.
      \label{eq:internal_equation}
\end{equation}
Whenever coordinate $i$ is selected, the post-update residual satisfies $r_i=0$. Under coordinate cyclicity, at most $T$ residual differences separate any sufficiently late iterate from the latest update of coordinate $i$. Equation \eqref{eq:internal_equation} therefore implies $r_i^{(k)} \to 0$ for every $i$. Hence $\mathbf{r}^{(k)} \to \mathbf{0}_{n \times 1}$ and, since $\mathbf{A}$ is non-singular, $\mathbf{x}^{(k)}\to \mathbf{A}^{-1} \mathbf{b} = \mathbf{x}^\star$.$\hfill \blacksquare$

\begin{proposition}
At step $k$, after the exact scaled update associated with coordinate $i$, we have that: i) $\mathbf{e}_i^T \mathbf{r}^{(k+1)} = 0$ and $(\mathbf{x}^{(k)})^T \mathbf{r}^{(k+1)} = 0$, and consequently ii) $( \mathbf{x}^{(k+1)})^T \mathbf{r}^{(k+1)} = 0$, i.e., the new iterate is also orthogonal to the new residual, and finally iii) $(\mathbf{r}^{(k+1)})^T \mathbf{A}^{-1}(\mathbf{r}^{(k+1)} - \mathbf{r}^{(k)}) = 0 $, i.e., for two consecutive residuals we have orthogonality in the $\mathbf{A}^{-1}$-inner product between the new residual and the residual difference. Furthermore, let $\mathbf{x}^\star=\mathbf{A}^{-1}\mathbf{b}$ and considering that $\mathbf{b}^T \mathbf{x}^{(k)} = (\mathbf{x}^{(k)})^T \mathbf{Ax}^{(k)}$, then the next decompositions hold:
\begin{equation}
\begin{aligned}
    \|\mathbf{x}^\star-\mathbf{x}^{(k)} \|_{\mathbf{A}}^2 & = \|\mathbf{x}^\star \|_{\mathbf{A}}^2 -\| \mathbf{x}^{(k)} \|_{\mathbf{A}}^2,\\
    \| \mathbf{x}^\star - \mathbf{x}^{(k+1)}  \|_{\mathbf{A}}^2 & = \| \mathbf{x}^\star  - \mathbf{x}^{(k)} \|_{\mathbf{A}}^2 - \| \mathbf{x}^{(k+1)} - \mathbf{x}^{(k)}  \|_{\mathbf{A}}^2, \\
    \| 
    \mathbf{r}^{(k+1)} \|_{\mathbf{A}^{-1}}^2 & = \| 
    \mathbf{r}^{(k)} \|_{\mathbf{A}^{-1}}^2 - \| \mathbf{r}^{(k+1)} - \mathbf{r}^{(k)}\|_{\mathbf{A}^{-1}}^2.
\end{aligned}
\end{equation}
\end{proposition}
\noindent\textbf{Proof.}
The exact update minimizes the quadratic over the subspace generated by $\mathbf{x}^{(k)}$ and $\mathbf{e}_i$. The gradient at the new iterate is therefore orthogonal to both search directions giving $i)$ and $ii)$. We also note that $\mathbf{x}^{(k+1)} - \mathbf{x}^{(k)} \in \operatorname{span}\{\mathbf{x}^{(k)},\mathbf{e}_i\}$ and therefore $(\mathbf{x}^{(k+1)} - \mathbf{x}^{(k)})^T \mathbf{r}^{(k+1)} = 0$ which together with $\mathbf{r}^{(k+1)} - \mathbf{r}^{(k)} = \mathbf{A}(\mathbf{x}^{(k+1)} - \mathbf{x}^{(k)})$ leads to $iii)$. The first decomposition follows from the expanded squared $\mathbf{A}$-norm and uses the scaling invariant \eqref{eq:scaling_invariant}. The last two decompositions follow because the new iterate is a linear combination of the two subspace directions. \hfill$\blacksquare$

\subsection{Geometric interpretation and numerical stability}

The numerator in~\eqref{eq:new_objective} is the squared directional derivative along $\mathbf{e}_i$. The denominator:
\begin{equation}
    A_{ii}-
    \frac{(\mathbf{a}_i^T\mathbf{x})^2}
    {\mathbf{x}^T\mathbf{Ax}},
    \label{eq:denominator}
\end{equation}
is the curvature along the component of $\mathbf{e}_i$ that is $\mathbf{A}$-orthogonal to $\mathbf{x}$. It depends on the current iterate and is no larger than the fixed diagonal curvature $A_{ii}$ used by ordinary coordinate descent. 
% Indeed,
% \begin{equation}
%     \frac{(\mathbf{a}_i^T\mathbf{x}-b_i)^2}{A_{ii}}
%     \leq
%     \frac{(\mathbf{a}_i^T\mathbf{x}-b_i)^2}
%     {A_{ii}-\dfrac{(\mathbf{a}_i^T\mathbf{x})^2}
%     {\mathbf{x}^T\mathbf{Ax}}}.
% \end{equation}
Furthermore, we have that:
\begin{equation}
    \gamma_i=\cos^2\varphi_i
    =\frac{(\mathbf{a}_i^T\mathbf{x})^2}
    {A_{ii}(\mathbf{x}^T\mathbf{Ax})},
\end{equation}
where $\varphi_i$ is the angle between $\mathbf{x}$ and $\mathbf{e}_i$ in the $\mathbf{A}$-inner product. Thus, the denominator \eqref{eq:denominator} can be written as $A_{ii}\sin^2\varphi_i$. Unfortunately, the quantity:
\begin{equation}
    A_{ii}(\mathbf{x}^T\mathbf{Ax})
    -(\mathbf{a}_i^T\mathbf{x})^2,
    \label{eq:local_gram_determinant}
\end{equation}
may be close to zero when $\mathbf{x}$ and $\mathbf{e}_i$ are nearly linearly dependent in the $\mathbf{A}$-inner product. In exact arithmetic, the quadratic objective remains bounded below. The numerical difficulty is instead that $\alpha^\star$ and $\beta^\star$ may become large and may cancel in the vector $\beta^\star\mathbf{x}+\alpha^\star\mathbf{e}_i$.

For non-zero $\mathbf{x}$, the denominator is zero when $\frac{(\mathbf{a}_i^T\mathbf{x})^2}
    {\mathbf{x}^T\mathbf{Ax}}
    =A_{ii},
    \text{ and therefore }
    \mathbf{x}=c\mathbf{A}^{-1}\mathbf{a}_i
    =c\mathbf{e}_i$. If~\eqref{eq:scaling_invariant} also holds, then the selected residual is zero, and the degenerate update can be skipped. Generally, the effective curvature can be written as:
\begin{equation}
    A_{ii}-
    \frac{(\mathbf{a}_i^T\mathbf{x})^2}
    {\mathbf{x}^T\mathbf{Ax}}
    =\mathbf{q}_i^T\mathbf{Aq}_i,
    \text{ where }
    \mathbf{q}_i
    =\mathbf{e}_i-
    \frac{\mathbf{a}_i^T\mathbf{x}}
    {\mathbf{x}^T\mathbf{Ax}}\mathbf{x}.
\end{equation}
Consequently we have $\mathbf{q}_i^T\mathbf{Aq}_i \geq\lambda_m(\mathbf{A})\|\mathbf{q}_i\|_2^2$, where $\|\mathbf{q}_i\|_2^2
    =1-2x_i
    \frac{\mathbf{a}_i^T\mathbf{x}}
    {\mathbf{x}^T\mathbf{Ax}}
    +\|\mathbf{x}\|_2^2
    \left(
    \frac{\mathbf{a}_i^T\mathbf{x}}
    {\mathbf{x}^T\mathbf{Ax}}
    \right)^2$.
This norm tends to zero as $\mathbf{x}$ approaches a multiple of $\mathbf{e}_i$.

% \noindent\textbf{Remark.}
% Let $\|\mathbf{x}\|_2=1$. For a fixed index $i$, the equality $\mathbf{a}_i^T\mathbf{x}=0$ is possible if and only if
% \begin{equation}
%     x_i^2\leq1-
%     \frac{A_{ii}^2}{\|\mathbf{a}_i\|_2^2}.
% \end{equation}

% \noindent\textbf{Proof.}
% Write $\mathbf{a}_i^T\mathbf{x}=A_{ii}x_i+\widetilde{\mathbf{a}}_i^T\mathbf{x}_{-i}$. If this quantity is zero, then
% \begin{equation}
%     |A_{ii}x_i|
%     \leq\|\widetilde{\mathbf{a}}_i\|_2
%     \sqrt{1-x_i^2}.
% \end{equation}
% Squaring and rearranging gives the stated inequality. Conversely, when the inequality holds, a vector $\mathbf{x}_{-i}$ with the required norm and inner product can be chosen. \hfill$\blacksquare$

% The ratio $A_{ii}^2/\|\mathbf{a}_i\|_2^2$ is the squared cosine of the Euclidean angle between $\mathbf{a}_i$ and $\mathbf{e}_i$ and is also related to the diagonal dominance of the $i$th column.

\begin{proposition}[Behavior of scaling near singularity]
Assume $ \mathbf{x}=c\mathbf{e}_i+\mathbf{z},
    \|\mathbf{x}\|_2=1,
    \|\mathbf{z}\|_2=\epsilon,
    \mathbf{e}_i^T\mathbf{z}=0$, with $0\leq\epsilon\leq1/\sqrt{\kappa(\mathbf{A})}$. Then we have:
\begin{equation}
    0\leq
    A_{ii}\frac{4\epsilon^2}{(\kappa(\mathbf{A})+1)^2}
    \leq
    A_{ii}-
    \frac{(\mathbf{a}_i^T\mathbf{x})^2}
    {\mathbf{x}^T\mathbf{Ax}}
    \leq
    A_{ii}\kappa(\mathbf{A})\epsilon^2
    \leq A_{ii}.
    \label{eq:near_coordinate_bounds}
\end{equation}
\end{proposition}
\noindent\textbf{Proof.}
Without loss of generality, isolate row and column $i$: $\mathbf{A}=
    \begin{bmatrix}
        A_{ii}&\widetilde{\mathbf{a}}_i^T\\
        \widetilde{\mathbf{a}}_i&\widetilde{\mathbf{A}}_i
    \end{bmatrix}$. Then $\mathbf{A}-\frac{1}{A_{ii}}\mathbf{a}_i\mathbf{a}_i^T
    =
    \begin{bmatrix}
        0&\mathbf{0}^T\\
        \mathbf{0}&\mathbf{S}
    \end{bmatrix}$ and $\mathbf{S}=\widetilde{\mathbf{A}}_i-
    \frac{1}{A_{ii}} \widetilde{\mathbf{a}}_i\widetilde{\mathbf{a}}_i^T$. The matrix $\mathbf{S}$ is the Schur complement of $A_{ii}$ and is positive definite. Its smallest eigenvalue satisfies
\begin{equation}
    \lambda_m(\mathbf{S})
    \geq
    \lambda_m(\widetilde{\mathbf{A}}_i)
    \left(
    1-
    \frac{\widetilde{\mathbf{a}}_i^T
     \widetilde{\mathbf{A}}_i^{-1}
    \widetilde{\mathbf{a}}_i}{A_{ii}}
    \right) =\frac{\lambda_m(\widetilde{\mathbf{A}}_i)}
    {A_{ii}(\mathbf{A}^{-1})_{ii}}
    \geq
    \frac{\lambda_m(\mathbf{A})}
    {A_{ii}(\mathbf{A}^{-1})_{ii}}.
\end{equation}
Moreover $A_{ii}-
    \frac{(\mathbf{a}_i^T\mathbf{x})^2}
    {\mathbf{x}^T\mathbf{Ax}}
    =A_{ii}
    \frac{\mathbf{z}^T\mathbf{Sz}}
    {\mathbf{x}^T\mathbf{Ax}}$. The upper estimate follows from
$\lambda_M(\mathbf{S}) \leq \lambda_M(\mathbf{A})$. For the lower estimate, use the preceding bound on $\lambda_m(\mathbf{S})$ and the Kantorovich inequality $A_{ii}(\mathbf{A}^{-1})_{ii}
    \leq \frac{(\kappa(\mathbf{A}) + 1)^2}
    {4\kappa(\mathbf{A})}.\hfill\blacksquare$
%This gives~\eqref{eq:near_coordinate_bounds}. \hfill$\blacksquare$

Practically, in Algorithm \ref{alg:scaled_coordinate_descent}, when \eqref{eq:local_gram_determinant} is too small, the safeguarded branch of Algorithm~1 first uses the ordinary coordinate correction and then rescales the entire corrected vector.

\subsection{The effects of the scaling factor}

We would like to quantify the effect of the scaling factor through the quantity
$\gamma_i=\frac{(\mathbf{a}_i^T\mathbf{x})^2}
{A_{ii}(\mathbf{x}^T\mathbf{A}\mathbf{x})}$. First, note that $\gamma_i$
is invariant under any nonzero scaling of $\mathbf{x}$. Therefore, the scaling component of the proposed update does not change $\gamma_i$ by itself; changes in $\gamma_i$ are caused by the coordinate correction, which changes the direction of the iterate.

In the nondegenerate case in which $\mathbf{x}$ is not proportional to $\mathbf{e}_i$, the Cauchy--Schwarz inequality in the $\mathbf{A}$-inner product gives $0\leq\gamma_i<1$. Moreover, $\gamma_i=0$ if and only if $\mathbf{a}_i^T \mathbf{x}=0$. Thus, $\gamma_i$ can vanish whenever $\mathbf{x}$ belongs to the $(n-1)$-dimensional hyperplane orthogonal to $\mathbf{a}_i$, and no strictly positive lower bound on $\gamma_i$ can hold uniformly for all iterates.

% Let $r_i=\mathbf{a}_i^T\mathbf{x}-b_i$.
The objective decrease produced by the exact scaled coordinate update is $\Delta_{\mathrm{SCD}}(i)
    =
    \frac{r_i^2}{2A_{ii}(1-\gamma_i)}$, whereas the decrease produced by ordinary coordinate descent is $\Delta_{\mathrm{CD}}(i)
    =
    \frac{r_i^2}{2A_{ii}}$.
Consequently, whenever $r_i\neq0$,
\begin{equation}
    \frac{\Delta_{\mathrm{SCD}}(i)}
    {\Delta_{\mathrm{CD}}(i)}
    =
    \frac{1}{1-\gamma_i},
    \ 
    \frac{\Delta_{\mathrm{SCD}}(i)
    -\Delta_{\mathrm{CD}}(i)}
    {\Delta_{\mathrm{CD}}(i)}
    =
    \frac{\gamma_i}{1-\gamma_i}.
\end{equation}
Therefore, the scaled method gives exactly the same one-step decrease as
ordinary coordinate descent when $\gamma_i=0$, and it gives a strictly
larger decrease whenever $r_i\neq0$ and $\gamma_i>0$. When
$\gamma_i\approx0$, the two methods have approximately the same
one-step behavior.

Since $\mathbf{a}_i^T\mathbf{x}=r_i+b_i$, the triangle and reverse
triangle inequalities give
\begin{equation}
    \bigl||r_i|-|b_i|\bigr|
    \leq
    |\mathbf{a}_i^T\mathbf{x}|
    \leq
    |r_i|+|b_i| \text{ it follows that }
    \gamma_i
    \geq
    \frac{\bigl(|r_i|-|b_i|\bigr)^2}
    {A_{ii}(\mathbf{x}^T\mathbf{A}\mathbf{x})}.
\end{equation}
Hence, if $r_i\neq0$ and $|r_i|\neq|b_i|$, then
$\mathbf{a}_i^T\mathbf{x}\neq0$, so $\gamma_i>0$ and the scaled update
gives a strictly larger objective decrease than ordinary coordinate
descent for the same selected coordinate. This condition is sufficient
but not necessary: even when $|r_i|=|b_i|$, it is possible to have
$\mathbf{a}_i^T\mathbf{x}\neq0$. Therefore, the exact benefit of scaling
is governed by $\gamma_i$, while the inequality above provides only a
state-dependent sufficient lower bound. We also highlight the following worst-case situation for our approach.

\begin{proposition}[Collapse to ordinary coordinate descent] Assume that the current iterate $\mathbf{x}^{(k)} \neq \mathbf{0}_{n \times 1}$ is optimally scaled, i.e., $(\mathbf{x}^{(k)})^T \mathbf{r}^{(k)} = 0$, and that, for the selected coordinate $i$, $\mathbf{a}_i^T\mathbf{x}^{(k)} = 0$. Then, the following hold:
\begin{enumerate}
    \item[i)] For $b_i\neq 0$, the scaled coordinate update reduces to the ordinary coordinate update, with $\beta^\star=1$ and $\alpha^\star = b_i/A_{ii}$. Moreover, $\mathbf{A}\mathbf{x}^{(k+1)} = b_i \mathbf{e}_i$ if and only if $\mathbf{x}^{(k)} = b_i\left( \mathbf{A}^{-1} \mathbf{e}_i - A_{ii}^{-1} \mathbf{e}_i \right)$. This candidate iterate satisfies the optimal scaling condition if and only if $\mathbf{b}^T\mathbf{A}^{-1}\mathbf{e}_i=b_i(\mathbf{A}^{-1})_{ii}$.

    \item[ii)] The compatibility conditions $\mathbf{b}^T\mathbf{A}^{-1}\mathbf{e}_i=b_i(\mathbf{A}^{-1})_{ii}$ hold for all coordinates $i$ if and only if $\operatorname{off}(\mathbf{A}^{-1})\mathbf{b}=\mathbf{0}_{n\times 1}$. Equivalently, $\mathbf{A}^{-1}\mathbf{b}=\operatorname{diag}(\operatorname{diag}(\mathbf{A}^{-1}))\mathbf{b}$. Thus, for such right-hand sides, the diagonal part of $\mathbf{A}^{-1}$ reproduces the action of the full inverse exactly. Non-zero right-hand sides satisfying this condition exist if and only if the hollow matrix $\operatorname{off}(\mathbf{A}^{-1})$ is singular. A trivial case is when $\mathbf{A}$ is diagonal.

    \item[iii)] Assume that $\mathbf{x}^{(1)} = b_1\mathbf{A}^{-1} \mathbf{e}_1 \neq \mathbf{0}_{n \times 1}$, so that $\mathbf{A} \mathbf{x}^{(1)} = b_1 \mathbf{e}_1$, and that $\mathbf{x}^{(1)}$ satisfies the optimal scaling condition. Suppose that the coordinates are subsequently selected in the order $2,3,\ldots,n$ and that, at each of the steps $2,\ldots,t-1$, the corresponding correlation is zero, so that the scaled update reduces to the ordinary coordinate update. Then, the same collapse occurs at coordinate $t$ if and only if $\sum_{j=2}^{t-1}\frac{A_{tj}}{A_{jj}}b_j=0$, for $t=3,\ldots,n$.

    \item[iv)] For the fixed coordinate ordering $1,2,\ldots,n$, the relations $\sum_{j=2}^{t-1}A_{tj}b_j/A_{jj}=0$, for $t=3,\ldots,n$, hold for every choice of $b_2,\ldots,b_{n-1}$ if and only if the principal submatrix of $\mathbf{A}$ indexed by $\{2,\ldots,n\}$ is diagonal. If $n\geq 3$ and the analogous property holds for every permutation of the coordinate ordering, then it holds if and only if $\mathbf{A}$ is diagonal.
\end{enumerate}
\end{proposition}
\noindent \textbf{Proof.} Since $\mathbf{a}_i^T\mathbf{x}^{(k)}=0$ and $(\mathbf{x}^{(k)})^T \mathbf{A} \mathbf{x}^{(k)} = \mathbf{b}^T \mathbf{x}^{(k)}$, the expressions for the optimal parameters reduce to $\beta^\star = 1$ and $\alpha^\star = b_i/A_{ii}$. Consequently, $\mathbf{x}^{(k+1)} = \mathbf{x}^{(k)}+(b_i/A_{ii})\mathbf{e}_i$. Therefore, $\mathbf{A}\mathbf{x}^{(k+1)} = b_i\mathbf{e}_i$ if and only if $\mathbf{A} \mathbf{x}^{(k)} = b_i \mathbf{e}_i -(b_i/A_{ii}) \mathbf{a}_i$, which, after multiplication by $\mathbf{A}^{-1}$, is equivalent to $\mathbf{x}^{(k)} = b_i\left(\mathbf{A}^{-1} \mathbf{e}_i - A_{ii}^{-1} \mathbf{e}_i \right)$. For this candidate iterate, $(\mathbf{x}^{(k)})^T\mathbf{A}\mathbf{x}_k=b_i^2((\mathbf{A}^{-1})_{ii}-1/A_{ii})$, whereas $\mathbf{b}^T \mathbf{x}^{(k)} = b_i (\mathbf{b}^T \mathbf{A}^{-1} \mathbf{e}_i - b_i/A_{ii})$. Since $b_i\neq 0$, the equality $(\mathbf{x}^{(k)})^T \mathbf{A} \mathbf{x}^{(k)} = \mathbf{b}^T \mathbf{x}^{(k)}$ holds if and only if $\mathbf{b}^T \mathbf{A}^{-1} \mathbf{e}_i = b_i(\mathbf{A}^{-1})_{ii}$. This proves i). Since $\mathbf{A}^{-1}$ is symmetric, $\mathbf{b}^T \mathbf{A}^{-1} \mathbf{e}_i = \mathbf{e}_i^T \mathbf{A}^{-1} \mathbf{b} = (\mathbf{A}^{-1} \mathbf{b})_i$. Hence, imposing the compatibility condition for every $i$ gives $\mathbf{A}^{-1} \mathbf{b} = \operatorname{diag}(\operatorname{diag}(\mathbf{A}^{-1})) \mathbf{b}$, or, equivalently, $\operatorname{off}(\mathbf{A}^{-1}) \mathbf{b} = \mathbf{0}_{n\times 1}$. A nonzero vector $\mathbf{b}$ satisfying this equality exists if and only if $\operatorname{off}(\mathbf{A}^{-1})$ has a nontrivial nullspace. This proves ii). For iii), assume that the zero-correlation collapse has occurred successively at coordinates $2,\ldots,t-1$. At each such step, $\beta^\star = 1$ and $\alpha^\star = b_j/A_{jj}$. Therefore, immediately before coordinate $t$ is selected at iteration $k$, $\mathbf{x}^{(k)} = b_1 \mathbf{A}^{-1} \mathbf{e}_1 + \sum_{j=2}^{t-1} b_j / A_{jj} \mathbf{e}_j$. Multiplication by $\mathbf{A}$ gives $\mathbf{A}\mathbf{x}^{(k)} = b_1 \mathbf{e}_1 + \sum_{j=2}^{t-1}(b_j/A_{jj}) \mathbf{a}_j$. Consequently, $\mathbf{a}_t^T \mathbf{x}^{(k)} = \mathbf{e}_t^T \mathbf{A} \mathbf{x}^{(k)} = \sum_{j=2}^{t-1} A_{tj}b_j/A_{jj}$. Thus, the zero-correlation condition $\mathbf{a}_t^T \mathbf{x}^{(k)} = 0$ holds if and only if $\sum_{j=2}^{t-1}A_{tj}b_j/A_{jj} = 0$. This proves iii). Finally, suppose that the relations in iii) hold for every choice of $b_2,\ldots,b_{n-1}$. For $t=3$, the equality $A_{32}b_2/A_{22}=0$ for every $b_2$ implies $A_{32}=0$. For $t=4$, the equality $A_{42}b_2/A_{22}+A_{43}b_3/A_{33}=0$ for every $b_2$ and $b_3$ implies $A_{42}=A_{43}=0$. Continuing in this way yields $A_{tj}=0$ for all $2\leq j<t\leq n$. By symmetry, the principal submatrix indexed by $\{2,\ldots,n\}$ is diagonal. The converse follows immediately.

If the analogous property holds for every permutation of the coordinate ordering, then any two distinct coordinates can be placed in the second and third positions. The corresponding case $t=3$ forces their off-diagonal entry to be zero. Hence, every off-diagonal entry of $\mathbf{A}$ is zero, and $\mathbf{A}$ is diagonal. The converse is immediate. $\hfill\blacksquare$

The proposition characterizes structured situations in which the additional
scaling direction becomes locally ineffective. At an individual iteration,
the condition $\mathbf{a}_i^T\mathbf{x}^{(k)}=0$ forces the scaled update to
coincide with the ordinary coordinate update. Requiring the associated
compatibility relation to hold for every coordinate is substantially more
restrictive: it identifies the exceptional subspace of right-hand sides
satisfying
$\operatorname{off}(\mathbf{A}^{-1})\mathbf{b} = \mathbf{0}_{n \times 1}$, on which the
diagonal part of $\mathbf{A}^{-1}$ reproduces the action of the complete
inverse. The cyclic conditions in $iii)$--$iv)$ further show that a persistent
collapse of scaled coordinate descent to ordinary coordinate descent is not
a generic phenomenon; it requires special cancellations depending on
$\mathbf{A}$, $\mathbf{b}$, and the coordinate ordering. If the collapse is
required for every right-hand side and every coordinate ordering, then
$\mathbf{A}$ must be diagonal. Thus, except for these structurally
degenerate cases, the scaling direction generally contributes an additional
degree of freedom and can produce a strictly larger one-step objective
decrease than ordinary coordinate descent.
%This result describes cases where scaling does not improve the ordinary coordinate descent, beyond the observation that $\mathbf{a}_i^T \mathbf{x}^{(k)}  = 0$ does not bring any benefits. Also, it provides a way to solve in $n$ step the simplest case, when $\mathbf{A}$ is a diagonal matrix: start from $\mathbf{x}^{(1)} = \mathbf{0}_{n \times 1}$, then $\mathbf{Ax}^{(2)} = b_1 \mathbf{e}_1$ by the ordinary coordinate descent step and then all the following steps $t \geq 3$ also perform an ordinary coordinate descent. Then this trivial problem is solved in exactly $n$ steps.

\subsection{Randomized selection and expected convergence}
\label{sec:expectation}

Given the current solution $\mathbf{x}$, define the probability of choosing index $i$:
\begin{equation}
    \mathbb{P}(I = i \mid \mathbf{x})=
    \frac{A_{ii} -
    \dfrac{(\mathbf{a}_i^T \mathbf{x})^2}
    {\mathbf{x}^T\mathbf{Ax}}}
    {\operatorname{tr}(\mathbf{A}) -
    \dfrac{\mathbf{x}^T \mathbf{A}^2 \mathbf{x}}
    {\mathbf{x}^T \mathbf{Ax}}}, i=1,\ldots,n,
    \label{eq:prob_spd}
\end{equation}
and this is normalized by the quantity:
\begin{equation}
    \sum_{i=1}^n
    \left(
    A_{ii}-
    \frac{(\mathbf{a}_i^T\mathbf{x})^2}
    {\mathbf{x}^T\mathbf{Ax}}
    \right)
    =\operatorname{tr}(\mathbf{A})-
    \frac{\mathbf{x}^T\mathbf{A}^2\mathbf{x}}
    {\mathbf{x}^T\mathbf{Ax}}.
\end{equation}
We also define:
\begin{equation}
    \delta(\mathbf{x})
    =f(\mathbf{x})-f(\mathbf{x}^\star)
     =  \frac12 \| \mathbf{x}^\star - \mathbf{x} \|_{\mathbf{A}}^2 =\frac{1}{2}\mathbf{r}^T\mathbf{A}^{-1}\mathbf{r},
    \text{ where }\mathbf{r}=\mathbf{Ax}-\mathbf{b},
\end{equation}
then we have the following convergence result.

\begin{proposition}[Probabilistic convergence result]
Assume $n\geq2$ and $\mathbf{x}^{(k)} \neq \mathbf{0}_{n \times 1}$ such that it is optimally scaled, i.e., $(\mathbf{x}^{(k)})^T \mathbf{r}^{(k)} = 0$, and the sampling probabilities in~\eqref{eq:prob_spd} for the indices, then we have that:
\begin{equation}
\begin{aligned}
    \mathbb{E}
    \left[f(\mathbf{x}^{(k+1)})\mid\mathbf{x}^{(k)}\right]
    =f(\mathbf{x}^{(k)})
    -\frac{1}{2}
    \frac{\|\mathbf{A}\mathbf{x}^{(k)}-\mathbf{b}\|_2^2}
    {\operatorname{tr}(\mathbf{A})-
    \dfrac{({\mathbf{x}^{(k)}})^T\mathbf{A}^2\mathbf{x}^{(k)}}
    {({\mathbf{x}^{(k)}})^T\mathbf{A}\mathbf{x}^{(k)}}}.
\end{aligned}
\label{eq:adaptive_sampling_exact}
\end{equation}
Moreover, we have:
\begin{equation}
\begin{aligned}
    \mathbb{E}
    \left[\delta(\mathbf{x}^{(k+1)})\mid\mathbf{x}^{(k)}\right]
    &\leq
    \left(
    1-
    \frac{\lambda_m(\mathbf{A})
    ({\mathbf{x}^{(k)}})^T\mathbf{A}\mathbf{x}^{(k)}}
    {({\mathbf{x}^{(k)}})^T
    (\operatorname{tr}(\mathbf{A})\mathbf{A}-\mathbf{A}^2)
    \mathbf{x}^{(k)}}
    \right)
    \delta(\mathbf{x}^{(k)})\\
    &\leq
    \left(
    1-
    \frac{\lambda_m(\mathbf{A})}
    {\operatorname{tr}(\mathbf{A})-\lambda_m(\mathbf{A})}
    \right)
    \delta(\mathbf{x}^{(k)}).
\end{aligned}
\label{eq:adaptive_sampling}
\end{equation}
\end{proposition}
\noindent\textbf{Proof.}
For coordinate $i$, the exact decrease is $\frac{1}{2} r_i^2 \left( A_{ii}-\dfrac{(\mathbf{a}_i^T \mathbf{x})^2} {\mathbf{x}^T \mathbf{Ax}} \right)^{-1}$.
We plug in the expectation of selecting the $i^\text{th}$ coordinate \eqref{eq:prob_spd} and develop:
\begin{equation}
	\begin{aligned}
		\mathbb{E}[ f(\mathbf{x}^{(k+1)})  | \mathbf{x}^{(k)} ] & = f(\mathbf{x}^{(k)}) - \frac{1}{2}\sum_{i=1}^n  \mathbb{P}(I = i \ | \ \mathbf{x}^{(k)})  \frac{ ( \mathbf{a}_i^T \mathbf{x}^{(k)} - b_i )^2 }{ A_{ii} - \frac{ ( \mathbf{a}_i^T {\mathbf{x}^{(k)}} )^2 }{ 
        (\mathbf{x}^{(k)})^T \mathbf{A} \mathbf{x}^{(k)} } } \\
		& = f(\mathbf{x}^{(k)}) - \frac{1}{2} \frac{ \|  \mathbf{Ax}^{(k)} - \mathbf{b} \|_2^2 }{ \operatorname{tr}(\mathbf{A}) - \frac{ (\mathbf{x}^{(k)})^T \mathbf{A}^2 \mathbf{x}^{(k)} }{ (\mathbf{x}^{(k)})^T \mathbf{A} \mathbf{x}^{(k)} } } \\
		& \leq  f(\mathbf{x}^{(k)}) - \frac{ {\mathbf{x}^{(k)}}^T \mathbf{A} \mathbf{x}^{(k)} }{ \operatorname{tr}(\mathbf{A}) ( (\mathbf{x}^{(k)})^T \mathbf{A}\mathbf{x}^{(k)} ) - (\mathbf{x}^{(k)})^T \mathbf{A}^2 \mathbf{x}^{(k)} } \lambda_m(\mathbf{A}) \delta(\mathbf{x}^{(k)}) \\
		& = f(\mathbf{x}^{(k)}) - \frac{ (\mathbf{x}^{(k)})^T \mathbf{A} \mathbf{x}^{(k)} }{ (\mathbf{x}^{(k)})^T ( \operatorname{tr}(\mathbf{A})  \mathbf{A} -  \mathbf{A}^2) \mathbf{x}^{(k)} } \lambda_m(\mathbf{A}) \delta(\mathbf{x}^{(k)}) .
	\end{aligned}
\end{equation}
In the summation term we have assumed $A_{ii} - \frac{ ( \mathbf{a}_i^T {\mathbf{x}^{(k)}} )^2 }{ (\mathbf{x}^{(k)})^T \mathbf{A} \mathbf{x}^{(k)} } \neq 0$ since in this case the residual is also zero and therefore the coordinate is stationary and the probability of selecting it should be zero. We have also used $\frac{1}{2}\|  \mathbf{Ax}^{(k)} - \mathbf{b} \|_2^2\geq \lambda_m(\mathbf{A}) \delta(\mathbf{x}^{(k)}) $. This leads to the result in \eqref{eq:adaptive_sampling}. The key quantity we considered here is the spectrum of the matrix $\operatorname{tr}(\mathbf{A})\mathbf{I}_n - \mathbf{A}$. The eigenvalues of this matrix are $\operatorname{tr}(\mathbf{A}) - \lambda_i$. These quantities are strictly greater than zero for any symmetric positive definite $\mathbf{A}$. Therefore, we have the following:
\begin{equation}
	\operatorname{tr}(\mathbf{A}) - \lambda_M(\mathbf{A})  \leq \frac{ (\mathbf{x}^{(k)})^T ( \operatorname{tr}(\mathbf{A})  \mathbf{A} -  \mathbf{A}^2) \mathbf{x}^{(k)}  }{(\mathbf{x}^{(k)})^T \mathbf{A} \mathbf{x}^{(k)} } \leq \operatorname{tr}(\mathbf{A}) - \lambda_m(\mathbf{A}).% \leq (n-1) \kappa(\mathbf{A}) \lambda_m(\mathbf{A}) .
\end{equation}
This result holds because of the generalized Rayleigh quotient. $\hfill \blacksquare$
% Multiplying by~\eqref{eq:prob_spd} and summing over $i$ cancels the coordinate-dependent denominator and gives~\eqref{eq:adaptive_sampling_exact}. Furthermore, $\frac{1}{2}\| \mathbf{r} \|_2^2
%  \geq \lambda_m( \mathbf{A} ) \delta( \mathbf{x} )$. Finally, the generalized Rayleigh quotient gives:
% \begin{equation}
%     \operatorname{tr}(\mathbf{A})-\lambda_M(\mathbf{A})
%     \leq
%     \frac{\mathbf{x}^T
%     (\operatorname{tr}(\mathbf{A})\mathbf{A}-\mathbf{A}^2)
%     \mathbf{x}}
%     {\mathbf{x}^T\mathbf{Ax}} \leq
%     \operatorname{tr}(\mathbf{A})-\lambda_m(\mathbf{A}),
% \end{equation}
% which proves~\eqref{eq:adaptive_sampling}. \hfill$\blacksquare$

In order to sample with \eqref{eq:prob_spd}, maintaining $\mathbf{A}\mathbf{x}^{(k)}$ requires $O(n)$ additional storage. The product is updated as:
\begin{equation}
    \mathbf{A}\mathbf{x}^{(k+1)}
    =\beta^\star\mathbf{A}\mathbf{x}^{(k)}
    +\alpha^\star\mathbf{a}_i,
    \label{eq:updating_dot_product}
\end{equation}
and all probabilities can then be formed in $O(n)$ operations. An efficient implementation stores the iterate as $\mathbf{x}=\sigma\overline{\mathbf{x}}$ and caching $\mathbf{v}=\mathbf{A}\overline{\mathbf{x}}$, so that $\mathbf{A}\mathbf{x}=\sigma\mathbf{v}$. The dynamic weights are then $q_i=A_{ii}-(\sigma^2/\omega)v_i^2$. Forming all weights and their cumulative sum costs exactly $4n$ operations. After selecting coordinate $i$, compute $r_i=\sigma v_i-b_i$, $\alpha=-r_i/q_i$, $\sigma^+=\sigma-\alpha(\sigma^2/\omega)v_i$, and $t=\alpha/\sigma^+$. The new state is obtained from $y_i+t$ and $\mathbf{v}+t\mathbf{a}_i$, where the latter update costs $2n$ operations. Hence, the complete dynamically sampled iteration costs $6n+12$ scalar operations, excluding comparisons and random-number generation. The explicit iterate is recovered only when needed as $\mathbf{x}=\sigma\mathbf{y}$. Thus, for a dense implementation, the resulting iteration cost is approximately $6n$ scalar operations.% under the counting convention used here.% - computing \eqref{eq:updating_dot_product} takes $3n$ operations, but we avoid the $2n$ operations in Algorithm 1 step 6, as now we have available all the scalar products between $\mathbf{A}$ and $\mathbf{x}^{(k)}$.

\begin{remark}[On maintaining $\mathbf{Ax}^{(k)}$]
    While we need to keep track of $\mathbf{Ax}^{(k)}$ in order to sample according to \eqref{eq:prob_spd}, in the implementation we have noticed that having all these dot products available and avoiding the individual computation of $\mathbf{a}_i^T \mathbf{x}^{(k)}$ in Algorithm 1 line 6 can lead to significant running times gains due to block operations and caching.
\end{remark}

For ordinary randomized coordinate descent with probabilities proportional to $A_{ii}$, we note the well-known convergence result of~\cite{LeventhalLewis2010}:
\begin{equation}
    \mathbb{E}
    \left[\delta(\mathbf{x}^{(k+1)})
    \mid\mathbf{x}^{(k)}\right]
    \leq
    \left(
    1-\frac{\lambda_m(\mathbf{A})}{\operatorname{tr}(\mathbf{A})}
    \right)
    \delta(\mathbf{x}^{(k)}).
\end{equation}
To quantify the difference with this classic result, the ratio of the two worst-case contraction bounds after $k$ steps is approximately:
\begin{equation}
    \frac{
    \left(
    1-\dfrac{\lambda_m(\mathbf{A})}
    {\operatorname{tr}(\mathbf{A})-\lambda_m(\mathbf{A})}
    \right)^k}
    {
    \left(
    1-\dfrac{\lambda_m(\mathbf{A})}
    {\operatorname{tr}(\mathbf{A})}
    \right)^k}
    \approx
    \exp\left(
    -k\left(
    \frac{\lambda_m(\mathbf{A})}{\operatorname{tr}(\mathbf{A})}
    \right)^2
    \right),
\end{equation}
when $\lambda_m(\mathbf{A})/\operatorname{tr}(\mathbf{A})$ is small. From a computational perspective, the proposed probabilities depend on the current iterate $\mathbf{x}^{(k)}$ and therefore need to be computed at each step, while the probabilities of~\cite{LeventhalLewis2010} are static throughout the iterative process and need to be computed, in $O(n)$, only once.

\section{Scaled block coordinate descent}

The single-coordinate scaled update provides the most transparent setting in which to introduce the proposed scaling mechanism and develop its theoretical properties. From a practical perspective, however, updating one coordinate at a time may not fully exploit modern computational architectures, since block operations offer better data reuse, higher arithmetic intensity, and more efficient use of optimized matrix and linear-algebra kernels. It is therefore natural to extend the method by replacing the single coordinate direction with a block of coordinate directions and jointly optimizing the block coefficients together with the scaling factor. The resulting block scaled coordinate descent method preserves the main geometric and theoretical ideas of the single-coordinate formulation, while providing a more appropriate implementation in many applications.

\subsection{The scaled block approach}

We extend the scaling idea to a block of $d$ directions collected in a full column-rank matrix $\mathbf{U}\in\mathbb{R}^{n\times d}$. Ordinary block coordinate descent solves:
\begin{equation}
    \underset{\boldsymbol{\alpha} \in \mathbb{R}^d}
    {\text{minimize}}\ f(\mathbf{x} + \mathbf{U} \boldsymbol{\alpha}),
\end{equation}
with the optimal solution given by:
\begin{equation}
    \boldsymbol{\alpha}^\star
    =-(\mathbf{U}^T\mathbf{A}\mathbf{U})^{-1}
    \mathbf{U}^T(\mathbf{Ax}-\mathbf{b}),
    \label{eq:best_alpha_ddimensional_original}
\end{equation}
and the objective function value at the optimum is:
\begin{equation}
\begin{aligned}
    \min_{\boldsymbol{\alpha}}
    f(\mathbf{x}+\mathbf{U}\boldsymbol{\alpha})
    =f(\mathbf{x})
    -\frac{1}{2}
    (\mathbf{Ax}-\mathbf{b})^T
    \mathbf{U}
    (\mathbf{U}^T\mathbf{A}\mathbf{U})^{-1}
    \mathbf{U}^T
    (\mathbf{Ax}-\mathbf{b}).
\end{aligned}
\label{eq:ddimensional_search}
\end{equation}

When $\mathbf{U}_\mathcal{S}$ consists of the standard basis vectors indexed by a set $\mathcal{S}$, denote by $\mathbf{A}_{\mathcal{S}, \mathcal{S}} = \mathbf{U}^T_\mathcal{S} \mathbf{A U}_\mathcal{S}$ the corresponding principal submatrix and by $\mathbf{r}_{\mathcal{S}} = \mathbf{U}^T_\mathcal{S} \mathbf{r} =(\mathbf{Ax}-\mathbf{b})_{\mathcal{S}}$ the selected residual. Then:
\begin{equation}
    \boldsymbol{\alpha}^\star
    =-\mathbf{A}_{\mathcal{S}, \mathcal{S}}^{-1}\mathbf{r}_{\mathcal{S}},
\end{equation}
and the objective function value simplifies to:
\begin{equation}
    \min_{\boldsymbol{\alpha}}
    f(\mathbf{x}+\mathbf{U}_\mathcal{S}\boldsymbol{\alpha})
    =f(\mathbf{x})-
    \frac{1}{2}
    \mathbf{r}_{\mathcal{S}}^T
    \mathbf{A}_{\mathcal{S}, \mathcal{S}}^{-1}
    \mathbf{r}_{\mathcal{S}}.
\end{equation}

Now, the proposed scaled block problem is given by:
\begin{equation}
    \underset{\boldsymbol{\alpha} \in \mathbb{R}^d,
    \,\beta \in \mathbb{R}}
    {\text{minimize}}\ f(\beta\mathbf{x}+\mathbf{U}_\mathcal{S}\boldsymbol{\alpha}).
    \label{eq:qp_d}
\end{equation}
Denoting $\mathbf{g}_{\mathcal{S}} = \mathbf{U}^T_\mathcal{S} \mathbf{Ax} = (\mathbf{Ax} )_\mathcal{S}$ we have the following optimality result.

\begin{proposition}[Optimal scaled block updates]
\label{prop:11}
Assume~\eqref{eq:scaling_invariant} and the non-zero $\mathbf{x} \notin \operatorname{span}(\mathbf{U}_\mathcal{S})$. Then we have:
\begin{equation}
    \min_{\boldsymbol{\alpha},\beta} f(\beta\mathbf{x} + \mathbf{U}_\mathcal{S} \boldsymbol{\alpha})
    = f(\mathbf{x})-
    \frac{1}{2}
    \mathbf{r}_{\mathcal{S}}^T
    \mathbf{A}_{\mathcal{S}, \mathcal{S}}^{-1}
    \mathbf{r}_{\mathcal{S}} -
    \frac{1}{2}
    \frac{
    (\mathbf{g}_{\mathcal{S}}^T
    \mathbf{A}_{\mathcal{S}, \mathcal{S}}^{-1}
    \mathbf{r}_{\mathcal{S}})^2}
    {\mathbf{x}^T \mathbf{Ax} -
    \mathbf{g}_{\mathcal{S}}^T
    \mathbf{A}_{\mathcal{S}, \mathcal{S}}^{-1} \mathbf{g}_{\mathcal{S}}}.
\end{equation}
The optimal parameters are given by:
\begin{equation}
    \beta^\star
    = 1 + 
    \frac{\mathbf{g}_{\mathcal{S}}^T
     \mathbf{A}_{\mathcal{S}, \mathcal{S}}^{-1}
    \mathbf{r}_{\mathcal{S}}}
    {\mathbf{x}^T\mathbf{Ax}-
    \mathbf{g}_{\mathcal{S}}^T
    \mathbf{A}_{\mathcal{S}, \mathcal{S}}^{-1}
    \mathbf{g}_{\mathcal{S}}},
    \boldsymbol{\alpha}^\star
    =-\mathbf{A}_{\mathcal{S},\mathcal{S}}^{-1}
    \left(
    \mathbf{r}_{\mathcal{S}}+
    \frac{\mathbf{g}_{\mathcal{S}}^T
    \mathbf{A}_{\mathcal{S}, \mathcal{S}}^{-1}
    \mathbf{r}_{\mathcal{S}}}
    {\mathbf{x}^T\mathbf{Ax}-
    \mathbf{g}_{\mathcal{S}}^T
    \mathbf{A}_{\mathcal{S}, \mathcal{S}}^{-1}
    \mathbf{g}_{\mathcal{S}}}
    \mathbf{g}_{\mathcal{S}}
    \right).
\label{eq:best_block_solution}
\end{equation}
\end{proposition}
\noindent\textbf{Proof.}
Search over the subspace generated by $\mathbf{W} = \begin{bmatrix} \mathbf{x} & \mathbf{U}_\mathcal{S} \end{bmatrix}$. The exact subspace decrease is $\frac{1}{2}
    (\mathbf{Ax} - \mathbf{b})^T
    \mathbf{W}
    (\mathbf{W}^T \mathbf{A} \mathbf{W})^{-1}
    \mathbf{W}^T
    (\mathbf{Ax} - \mathbf{b})$. Applying the Schur complement inverse formula to the matrix in the preceding display and using the invariant~\eqref{eq:scaling_invariant} gives the stated decrease and the parameters in~\eqref{eq:best_block_solution}. The denominator is positive when $\mathbf{x} \notin \operatorname{span}(\mathbf{U}_\mathcal{S})$, and vanishes exactly when
$\mathbf{x}\in\operatorname{span}(\mathbf{U}_\mathcal{S})$. \hfill$\blacksquare$

%As in the one-coordinate case, the exact block update preserves $\mathbf{b}^T \mathbf{x} = \mathbf{x}^T \mathbf{Ax}$.

\subsection{Block algorithm and computational cost}
\begin{algorithm}[!th]
\caption{Block Scaled Coordinate Descent (BSCD)}
\label{alg:block_scaled_coordinate_descent}
\begin{algorithmic}[1]
\Require Symmetric positive definite matrix $\mathbf{A}\in\mathbb{R}^{n\times n}$, $\mathbf{b}\in\mathbb{R}^{n}$, $\mathbf{b}\neq\mathbf{0}$;
\Require $\mathbf{x}^{(0)}\neq\mathbf{0}$ satisfying $\mathbf{b}^{T}\mathbf{x}^{(0)}\neq0$;
\Require Block size $d$, maximum number of iterations $K$, threshold $\kappa_{\text{proj}} \!\geq\! 1$;
\Ensure Approximate solution $\mathbf{x}^{(K)}$.

\State Compute quantity $\omega^{(0)} \gets \mathbf{b}^{T} \mathbf{x}^{(0)}$;
\State Compute optimal scaling $\displaystyle \beta^{(0)}\gets
\frac{\omega^{(0)}}
{(\mathbf{x}^{(0)})^{T}\mathbf{A}\mathbf{x}^{(0)}}$;
\State Apply optimal scaling $\mathbf{x}^{(0)}\gets\beta^{(0)}\mathbf{x}^{(0)}$ and update $\omega^{(0)}\gets\beta^{(0)}\omega^{(0)}$.

\For{$k=1,\ldots,K$}
    \State Choose an index set $\mathcal{S}\subseteq\{1,\ldots,n\}$ with $|\mathcal{S}|=d$;
    \State Compute scalar products $\mathbf{g}_{\mathcal{S}} \gets \mathbf{A}_{:,\mathcal{S}}^{T} \mathbf{x}^{(k-1)}$, residuals $\mathbf{r}_{\mathcal{S}} \gets
    \mathbf{g}_{\mathcal{S}} - \mathbf{b}_{\mathcal{S}}$;
    \State Solve $\mathbf{A}_{\mathcal{S},\mathcal{S}} \mathbf{y}
     = \mathbf{g}_{\mathcal{S}}$ and $\mathbf{A}_{\mathcal{S},\mathcal{S}} \mathbf{z}
     = \mathbf{r}_{\mathcal{S}}$;

    %\If{$\frac{\omega^{(k-1)}-\mathbf{g}_{\mathcal{S}}^{T} \mathbf{y}}{\omega^{(k-1)}} \geq \frac{4 \kappa_{\text{proj}}}{ ( \kappa_{\text{proj}} + 1 )^2 }$}
    \If{$\mathbf{g}_{\mathcal{S}}^{T} \mathbf{y} \leq \left( \frac{ \kappa_{\text{proj}} -1 }{ \kappa_{\text{proj}} + 1  }\right)^2 \omega^{(k-1)}$}
     
    \State Compute quantity $\displaystyle
    \theta^{(k)} \gets \frac{\mathbf{g}_{\mathcal{S}}^{T} \mathbf{z}}
    {\omega^{(k-1)}-\mathbf{g}_{\mathcal{S}}^{T} \mathbf{y}}$;
    \State Compute $\beta^{(k)} \gets 1+\theta^{(k)}$ and $\boldsymbol{\alpha}^{(k)} \gets 
    - \mathbf{z} - \theta^{(k)} \mathbf{y}$;
    \State Scale $\mathbf{x}^{(k)}\gets
    \beta^{(k)}\mathbf{x}^{(k-1)}$ and update solution $\mathbf{x}_{\mathcal{S}}^{(k)} \gets 
    \mathbf{x}_{\mathcal{S}}^{(k)}
     + \boldsymbol{\alpha}^{(k)}$;
    \State Update quantity $\omega^{(k)} \gets \beta^{(k)} \omega^{(k-1)} + \mathbf{b}_{\mathcal{S}}^{T} \boldsymbol{\alpha}^{(k)}$;
    \Else
        \State Compute ordinary block coordinate correction
        $\boldsymbol{\alpha}^{(k)}\gets-\mathbf{z}$;
        \State Compute the optimal scaling
        $ \beta^{(k)}
        \gets
        \frac{
        \omega^{(k-1)}
        +\mathbf{b}_{\mathcal{S}}^{T}\boldsymbol{\alpha}^{(k)}
        }{
        \omega^{(k-1)}
        +2\mathbf{g}_{\mathcal{S}}^{T}\boldsymbol{\alpha}^{(k)}
        +\mathbf{r}_\mathcal{S}^T \mathbf{z}
        }$;

    %     \State Scale $\mathbf{x}^{(k)}\gets
    % \beta^{(k)}\mathbf{x}^{(k-1)}$, and update solution $\mathbf{x}_{\mathcal{S}}^{(k)} \gets 
    % \mathbf{x}_{\mathcal{S}}^{(k)}
    %  + \boldsymbol{\alpha}^{(k)}$;
        \State Update solution $\mathbf{x}^{(k)}
        \gets
        \beta^{(k)}
        (
        \mathbf{x}^{(k-1)}
        +\mathbf{U}_{\mathcal{S}}\boldsymbol{\alpha}^{(k)}
        )$;
        
        \State Update quantity $\omega^{(k)}
        \gets
        \beta^{(k)}
        (
        \omega^{(k-1)}
        +\mathbf{b}_{\mathcal{S}}^{T}\boldsymbol{\alpha}^{(k)}
        )$;
    \EndIf
\EndFor

\State \textbf{return} $\mathbf{x}^{(K)}$.
\end{algorithmic}
\end{algorithm}

Based on the optimality result in Proposition~\ref{prop:11}, Algorithm~\ref{alg:block_scaled_coordinate_descent} extends scaled coordinate descent by replacing the single coordinate direction with a block of $d$ coordinate directions. The initial vector is first optimally rescaled so that $\mathbf{b}^{T}\mathbf{x}^{(0)}=(\mathbf{x}^{(0)})^{T}\mathbf{A}\mathbf{x}^{(0)}=\omega^{(0)} > 0$ from the beginning.

At iteration $k$, a coordinate block $\mathcal{S}$ is selected and the vectors $\mathbf{g}_{\mathcal{S}}$ and $\mathbf{r}_{\mathcal{S}}$ are formed. Since $\mathbf{A}$ is symmetric positive definite, the principal submatrix $\mathbf{A}_{\mathcal{S},\mathcal{S}}$ is also symmetric positive definite, and the two systems defining $\mathbf{y}$ and $\mathbf{z}$ can be solved using the same factorization. The quantities $\theta$, $\beta$, and $\boldsymbol{\alpha}$ are then the optimal coefficients for minimizing the quadratic objective over the enlarged subspace generated by the current iterate and the selected coordinate block. The update jointly rescales the complete iterate and corrects the coordinates in $\mathcal{S}$. After an exact step, the new residual is orthogonal to the previous iterate and to every selected coordinate direction; consequently, the residual components indexed by $\mathcal{S}$ vanish and the invariant $\mathbf{b}^{T}\mathbf{x} = \mathbf{x}^T \mathbf{A} \mathbf{x} = \omega$ is preserved. In the computation of $\theta$, the denominator $\omega -\mathbf{g}_{\mathcal{S}}^{T} \mathbf{y} \geq 0$ is the Schur complement associated with the added scaling direction and is strictly positive whenever $\mathbf{x}$ does not belong to the span of the selected coordinate block. In the safeguarded step we use the same ideas as described in Algorithm~\ref{alg:scaled_coordinate_descent}, we perform an optimal block coordinate descent step with the set $\mathcal{S}$ and then we optimally scale the solution by:
\begin{equation}
    \beta^\dagger
        \gets
        \frac{
        \omega  + \mathbf{b}_{\mathcal{S}}^{T} \boldsymbol{\alpha}
        }{
        \omega + 2\mathbf{g}_{\mathcal{S}}^{T} \boldsymbol{\alpha}
        +\mathbf{r}_\mathcal{S}^T \mathbf{z}.
        }
\end{equation}

Generalizing the scalar safeguard from Algorithm~\ref{alg:scaled_coordinate_descent}, the new safeguard is designed to measure whether the current scaling direction is nearly contained in the selected coordinate block subspace. If we define $\gamma_{\mathcal{S}}
    =
    \frac{
    \mathbf{g}_{\mathcal{S}}^{T}\mathbf{y}
    }{
    \omega
    }$ then the vector $\mathbf{U}_{\mathcal{S}}\mathbf{y}$ is the
$\mathbf{A}$-orthogonal projection of $\mathbf{x}$ onto the selected
coordinate block subspace. Hence, $\gamma_{\mathcal{S}}$ is the fraction
of the $\mathbf{A}$-energy of $\mathbf{x}$ contained in that subspace, whereas $1-\gamma_{\mathcal{S}}
    =
    \frac{
    \omega-\mathbf{g}_{\mathcal{S}}^{T}\mathbf{y}
    }{
    \omega
    }$ is the relative squared $\mathbf{A}$-distance of $\mathbf{x}$ from the
block. The exact scaled update is accepted only when
$\gamma_{\mathcal{S}}\leq\gamma_{\max}<1$. Equivalently, choosing $\gamma_{\max}
    =
    \left(
    \frac{
    \kappa_{\text{proj}}-1
    }{
    \kappa_{\text{proj}}+1
    }
    \right)^2$ is identical to requiring that the normalized augmented Gram matrix have a condition number at most $\kappa_{\text{proj}}$. This safeguard
does not estimate the intrinsic condition number of
$\mathbf{A}_{\mathcal{S},\mathcal{S}}$ which would be computationally expensive to estimate in general; it controls only the additional near-dependence introduced by the scaling direction. If the containment threshold is exceeded, the algorithm uses the ordinary exact block correction and then optimally scales the complete corrected iterate. Similarly to Algorithm~\ref{alg:scaled_coordinate_descent}, the safeguarding step has the same effect on the indices belonging to the set $\mathcal{S}$ of the new residual which is now $(\beta^\dagger - 1) \mathbf{b}_\mathcal{S}$, and therefore non-zero in general.%we use the conditioning check, and when necessary we switch to the ordinary block coordinate update followed by the separate optimal scaling.

For a dense implementation, the operation counts can be precisely stated. Computing the selected products
$\mathbf{g}_{\mathcal{S}}
=\mathbf{A}_{:,\mathcal{S}}^{T}\mathbf{x}$ and the selected residual requires $2nd$ flops. Using a standard Cholesky factorization, the factorization of $\mathbf{A}_{\mathcal{S},\mathcal{S}}$ costs
$\frac13 d^3 + \frac12 d^2 + \frac16 d$ flops, and each right-hand side requires $2d^2$ flops for the two triangular solves. Consequently, an ordinary block coordinate update costs:
\begin{equation}
    N_{\text{BCD}}(n,d)
    =
    2nd+\frac{d^3}{3}
    +\frac{5d^2}{2}
    +\frac{7d}{6},
\end{equation}
flops. The scaled block method uses the same $d\times d$ factorization, but solves for two right-hand sides and additionally forms the scaling coefficients, scales the complete iterate, and updates
$\omega$ by $\omega-\boldsymbol{\alpha}^{T} \mathbf{r}_{\mathcal{S}}$. The total cost is:
\begin{equation}
    N_{\text{BSCD}}(n,d)
    =
    2nd+n+\frac{d^3}{3}
    +\frac{9d^2}{2}
    +\frac{55d}{6}+3,
\end{equation}
flops. Hence, for the same block size, scaling adds $n+2d^2+8d+3$ flops. The dominant additional term is one length-$n$ scaling when $d\ll n$, while the extra $2d^2$ term is the cost of solving the second
system with the already available Cholesky factors. In particular, the
scaled method does not require the factorization of a
$(d+1)\times(d+1)$ matrix; the same principal-block factorization is
reused for both right-hand sides. Furthermore, for comparison with ordinary BCD of size $d+1$ and thus matching the degrees of freedom we have that $N_{\mathrm{BCD}}(n,d+1) - N_{\mathrm{BSCD}}(n,d) = n-d^2 - 2d+ 1$. Hence, when $n > d^2
+2d - 1$ scaled block of size $d$ is arithmetically cheaper than an ordinary block of size $d+1$, despite both using a search space of dimension $d+1$.
%For a dense implementation, an ordinary block update is dominated by $N_\text{ordinary} = 2dn+\frac{1}{3}d^3+O(d^2)$  operations, whereas the scaled block update is dominated by $N_\text{scaled} = 2dn+n+\frac{1}{3}(d+1)^3+O((d+1)^2)$. Thus, relative to the ordinary method approach of the same size, scaling adds approximately $n+O(d^2)$ operations. For $d=1$, this recovers the comparison between approximately $2n$ operations for ordinary coordinate descent and $3n$ operations for the scaled update. Relative to a block of size $d+1$, the scaled block of size $d$ saves one column update but uses one dense scaling, so the leading-order costs differ by $O(n)$.

\subsection{A randomized convergence result}

While Algorithm~\ref{alg:block_scaled_coordinate_descent} uses a safeguarded update for badly conditioned practical scenarios, the result in this section refers only to the newly scaled updates developed in this paper.

To quantify the performance of the proposed method, we adapt the analysis performed in the sketch-and-project framework developed by Gower and Richtárik in~\cite{GowerRichtarik2015}. The next result establishes the convergence rate of the proposed method and compares it with the ordinary block coordinate descent.

\begin{proposition}[Convergence of scaled block coordinate descent]
Let $\mathbf{A}\in\mathbb{R}^{n\times n}$ be symmetric positive definite and $d \geq 1$ the block size. Let $\mathcal{S}, |\mathcal{S}| = d$, be sampled independently at every iteration from a fixed
distribution over coordinate subsets, and let $\mathbf{U}_{\mathcal{S}}$ denote the corresponding coordinate selection
matrix. Define $\mathbf{P}_{\mathcal{S}}^\text{BCD}
    =
    \mathbf{A}^{1/2} \mathbf{U}_{\mathcal{S}}
    (
    \mathbf{U}_{\mathcal{S}}^{T}
    \mathbf{A}
    \mathbf{U}_{\mathcal{S}}
    )^{-1}
    \mathbf{U}_{\mathcal{S}}^{T}
    \mathbf{A}^{1/2}$,
and assume that
    $\overline{\mathbf{P}}^\text{BCD}
    \! = \!
    \mathbb{E}_{\mathcal{S}}
    \left[
    \mathbf{P}_{\mathcal{S}}^\text{BCD}
    \right]
    \! \! \succ \! \! \mathbf{0}$.
Ordinary randomized block coordinate descent satisfies:
\begin{equation}
    \mathbb{E}
    \left[
    \delta(\mathbf{x}^{(k+1)})
    \mid
    \mathbf{x}^{(k)}
    \right]
    \leq
    (
    1 - \lambda_{m}
    (
    \overline{\mathbf{P}}^\text{BCD}
    )
    )
    \delta(\mathbf{x}^{(k)}).
    \label{eq:regular_block_convergence}
\end{equation}

Assume now that $\mathbf{x}^{(k)}$ is optimally scaled and define $\mathbf{W}_{\mathcal{S}}(\mathbf{x}^{(k)}) = \begin{bmatrix}
        \mathbf{x}^{(k)} &
        \mathbf{U}_{\mathcal{S}}
    \end{bmatrix}$.
For the exact scaled block update, define the state-dependent projection matrix
$\mathbf{P}_{\mathcal{S}}^{\text{BSCD}}
    (\mathbf{x}^{(k)})
    =
    \mathbf{A}^{1/2}
    \mathbf{W}_{\mathcal{S}}(\mathbf{x}^{(k)})
    (
    \mathbf{W}_{\mathcal{S}}(\mathbf{x}^{(k)})^{T}
    \mathbf{A}
    \mathbf{W}_{\mathcal{S}}(\mathbf{x}^{(k)})
    )^{\dagger}
    \mathbf{W}_{\mathcal{S}}(\mathbf{x}^{(k)})^{T}
    \mathbf{A}^{1/2}$.
Then:
\begin{equation}
    \mathbb{E}
    \left[
    \delta(\mathbf{x}^{(k+1)})
    \mid
    \mathbf{x}^{(k)}
    \right]
    \leq
    (
    1- \lambda_{m}
    (
    \mathbb{E}_{\mathcal{S}}
    \left[
    \mathbf{P}_{\mathcal{S}}^{\text{BSCD}}
    (\mathbf{x}^{(k)})
    \mid
    \mathbf{x}^{(k)}
    \right]
    )
    )
    \delta(\mathbf{x}^{(k)}).
\label{eq:scaled_block_convergence}
\end{equation}

For every block $\mathcal{S}$, let $\mathbf{w}_{\mathcal{S}}(\mathbf{x})
    =
    (
    \mathbf{I}-\mathbf{P}_{\mathcal{S}}^\text{BCD}
    )
    \mathbf{A}^{1/2}\mathbf{x}$. If $\mathbf{w}_{\mathcal{S}}(\mathbf{x})\neq\mathbf{0}$, then:
\begin{equation}
\mathbf{P}_{\mathcal{S}}^{\text{BSCD}}(\mathbf{x})
    =
    \mathbf{P}_{\mathcal{S}}^\text{BCD}
    +
    \frac{
    \mathbf{w}_{\mathcal{S}}(\mathbf{x})
    \mathbf{w}_{\mathcal{S}}(\mathbf{x})^{T}
    }{
    \|\mathbf{w}_{\mathcal{S}}(\mathbf{x})\|_{2}^{2}
    }.
\label{eq:scaled_block_rank_one_improvement}
\end{equation}
If $\mathbf{w}_{\mathcal{S}}(\mathbf{x})=\mathbf{0}_{n \times 1}$, then
$\mathbf{P}_{\mathcal{S}}^{\text{BSCD}}(\mathbf{x}) = \mathbf{P}_{\mathcal{S}}^\text{BCD}$. Hence, scaling augments each ordinary block projector by at most one positive semidefinite rank-one projector. Finally, because: $\mathbf{P}_{\mathcal{S}}^{\text{BSCD}}
    (\mathbf{x}^{(k)})
    \succeq
    \mathbf{P}_{\mathcal{S}}^\text{BCD}$ we have that:
\begin{equation}
    \lammin (\mathbb{E}[ \mathbf{P}_{\mathcal{S}}^{\text{BSCD}}
    (\mathbf{x}^{(k)}) | \mathbf{x}^{(k)}])
    \geq 
    \lammin (\mathbb{E}[ \mathbf{P}_{\mathcal{S}}^\text{BCD}]).
\end{equation}
\end{proposition}
\noindent \textbf{Proof.}
For ordinary randomized block coordinate descent, take
$\mathbf{B}=\mathbf{A}$ and
$\mathbf{S}=\mathbf{U}_{\mathcal{S}}$ in the sketch-and-project framework~\cite{GowerRichtarik2015}. Their randomized Newton
result gives the exact block coordinate update, and their Theorem~4.6 gives:
\begin{equation}
    \mathbb{E}
    \left[
    \|\mathbf{x}^{\star} - \mathbf{x}^{(k)}\|_{\mathbf{A}}^{2}
    \right]
    \leq
    (
    1-
    \lambda_{m}(\overline{\mathbf{P}}^\text{BCD})
    )^{k}
    \|\mathbf{x}^{\star} - \mathbf{x}^{(0)}\|_{\mathbf{A}}^{2}.
\end{equation}
Since
$\delta(\mathbf{x}) = \frac{1}{2} \|\mathbf{x}^{\star} - \mathbf{x}\|_{\mathbf{A}}^{2}$,
this gives~\eqref{eq:regular_block_convergence}. For the scaled method, condition on the current iterate
$\mathbf{x}^{(k)}$. The exact scaled block step minimizes the quadratic
over $\operatorname{span}
    \{
    \mathbf{x}^{(k)},
    \mathbf{U}_{\mathcal{S}}
    \}$. Equivalently, the transformed error
$\mathbf{A}^{1/2}
(\mathbf{x}^{\star} - \mathbf{x}^{(k)})$
is projected by
$\mathbf{I}
-\mathbf{P}_{\mathcal{S}}^{\text{SCD}}
(\mathbf{x}^{(k)})$.
Since
$\mathbf{P}_{\mathcal{S}}^{\text{SCD}}
(\mathbf{x}^{(k)})$
is an orthogonal projector, we have that $\mathbb{E}
    \left[
    \|\mathbf{x}^{\star} - \mathbf{x}^{(k+1)}\|_{\mathbf{A}}^{2}
    \mid
    \mathbf{x}^{(k)}
    \right]
    =
    \|\mathbf{x}^{\star}  - \mathbf{x}^{(k)} \|_{\mathbf{A}}^{2}
    -
    (
    \mathbf{A}^{1/2}
    (\mathbf{x}^{\star} - \mathbf{x}^{(k)})
    )^{T}
    \mathbb{E}_{\mathcal{S}}
    \left[
    \mathbf{P}_{\mathcal{S}}^{\text{SCD}}
    (\mathbf{x}^{(k)})
    \mid
    \mathbf{x}^{(k)}
    \right]
    (
    \mathbf{A}^{1/2}
    (\mathbf{x}^{\star} - \mathbf{x}^{(k)})
    )$. Then, bounding the quadratic form by its smallest eigenvalue gives
\eqref{eq:scaled_block_convergence}. This is the same conditional
projection argument used in the proof of Theorem~4.6, but the expectation
is conditioned on $\mathbf{x}^{(k)}$ because the scaled search space is
state-dependent. Finally, $\operatorname{range} ( \mathbf{A}^{1/2}\mathbf{U}_{\mathcal{S}}
    )
    \subseteq
    \operatorname{range}
    (
    \mathbf{A}^{1/2}
    \mathbf{W}_{\mathcal{S}}(\mathbf{x}^{(k)})
    )$.
The orthogonal projectors onto these nested subspaces therefore satisfy $\mathbf{P}_{\mathcal{S}}^{\text{BSCD}}
    (\mathbf{x}^{(k)})
    \succeq
    \mathbf{P}_{\mathcal{S}}^\text{BCD}$. Taking expectations and using the monotonicity of the minimum eigenvalue proves the result. The rank-one identity
\eqref{eq:scaled_block_rank_one_improvement} follows by orthogonally
decomposing the newly added vector
$\mathbf{A}^{1/2}\mathbf{x}^{(k)}$ relative to
$\operatorname{range}
(\mathbf{A}^{1/2}\mathbf{U}_{\mathcal{S}})$.$\hfill\blacksquare$

Thus, under the same block selection distribution, scaled randomized block coordinate descent inherits the complete convergence guarantee of ordinary randomized block coordinate descent and has a state-dependent contraction factor that is never worse. 

\section{Experimental results}

In this section, we perform several experiments to compare the
convergence rate of the proposed method against ordinary coordinate descent approaches, using various coordinate selection methods.

We investigate whether the additional decrease obtained by scaling the current iterate compensates for its computational cost. The experiments use synthetic symmetric positive definite quadratic problems with $n=5000$, and the results are averaged over $100$ realizations of the problem data and initial vector. Within each realization, the compared methods solve the same problem and use a common starting vector, optimally rescaled. Accuracy is measured by the objective gap $\delta(\mathbf{x})=f(\mathbf{x})-f(\mathbf{x}^{\star})$, where $\mathbf{x}^{\star}=\mathbf{A}^{-1}\mathbf{b}$. We compare this gap against cumulative elapsed time, so that the cost of coordinate selection and the additional work of the scaled update are reflected in the comparison.

For a given size $n$, all test matrices $\mathbf{A}$ and target vectors $\mathbf{b}$ are generated randomly from the standard Gaussian distribution, and $\mathbf{A}$ is symmetrized by summation. In each case, we report progress by both iteration count and running times.

\begin{figure}[!h]
    \centering    \includegraphics[width=0.9\linewidth]{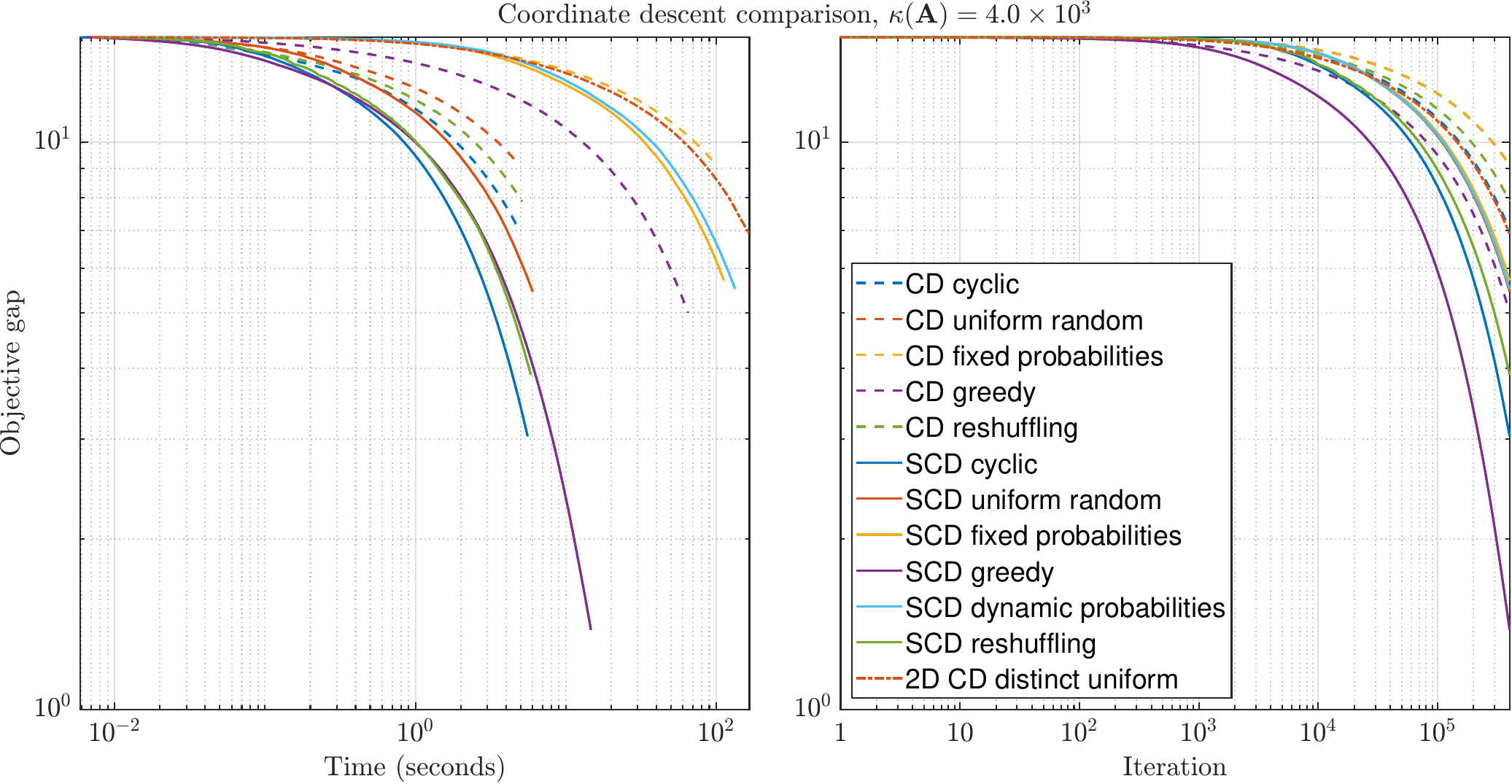}
    %\vspace{0.4cm}
    \includegraphics[width=0.9\linewidth]{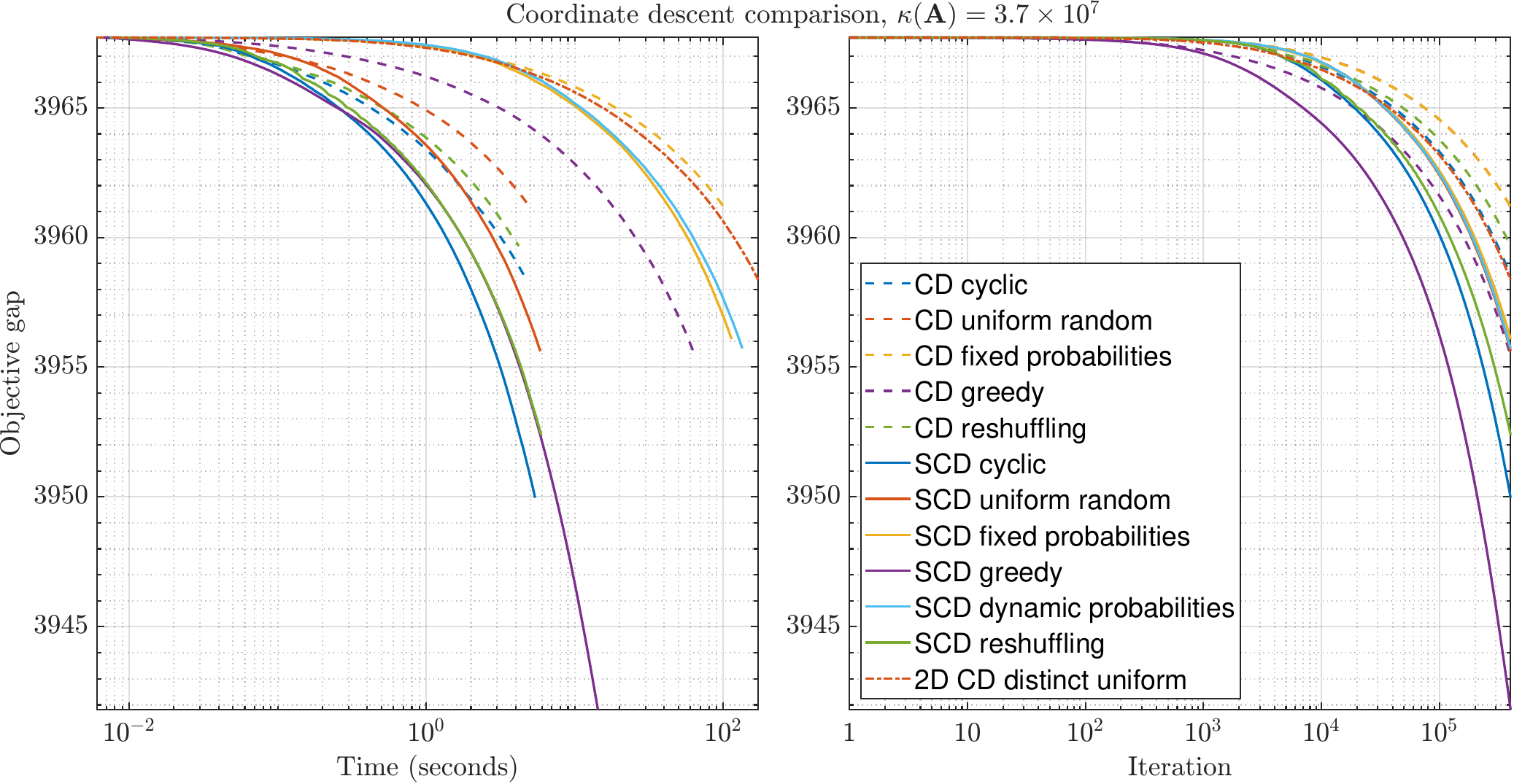}

    \caption{Comparison of the coordinate descent methods for two test problems with different condition numbers averaged over 100 realizations, for all $\mathbf{A}$, $\mathbf{b}$ and $\mathbf{x}^{(0)}$.}\label{fig:regular_n5000_comparison}
\end{figure}
\subsection{Single coordinate methods}

Figure~\ref{fig:regular_n5000_comparison} compares ordinary coordinate descent, the proposed scaled variants, and the two-coordinate baselines. The two configurations have reported average condition numbers $\kappa(\mathbf{A})\approx4.0134\times10^{3}$ in the upper plot and $\kappa(\mathbf{A})\approx3.6977\times10^{7}$ in the lower plot, with an iteration budget of $K=400\,000$. These configurations examine the effect of scaling in two substantially different conditioning regimes.

The coordinate selection rules have different computational requirements. The simplest index selection rules are \textit{cyclic} in order and \textit{reshuffling} in random order, ensuring all indices are visited at least once per sweep. \textit{Fixed probability} selection uses $p_i=A_{ii}/\operatorname{tr}(\mathbf{A})$, whereas \textit{dynamically sampled} SCD uses the state-dependent probabilities in~\eqref{eq:prob_spd}. \textit{Greedy selection} maximizes the current objective decrease, using the score $r_i^2/A_{ii}$ for ordinary CD and $r_i^2/(A_{ii}(1-\gamma_i))$ for the exact scaled update. The two-coordinate methods provide a complementary comparison: they enlarge the local search space using a second coordinate, whereas SCD uses the current iterate as its additional direction. Thus, the comparison concerns both the choice of search space and the cost of selecting it. The same coordinate selection rules are considered for both CD and the proposed SCD. We also show results from the \textit{two-dimensional ordinary coordinate descent} method, where the two distinct indices are selected uniformly at random and independently.

The exact one-step dominance established earlier applies to a common current iterate and a common selected coordinate. It does not imply that independently generated trajectories must remain ordered, or that SCD must be faster for every selection rule. The comparison in terms of the elapsed time is therefore essential: a larger decrease per iteration is useful computationally only when it compensates for the additional work. Cyclic, reshuffling, and greedy selection criteria for SCD establish the best empirical results in our simulations.

\begin{figure}[!h]
    \centering

    \includegraphics[width=0.9\linewidth]{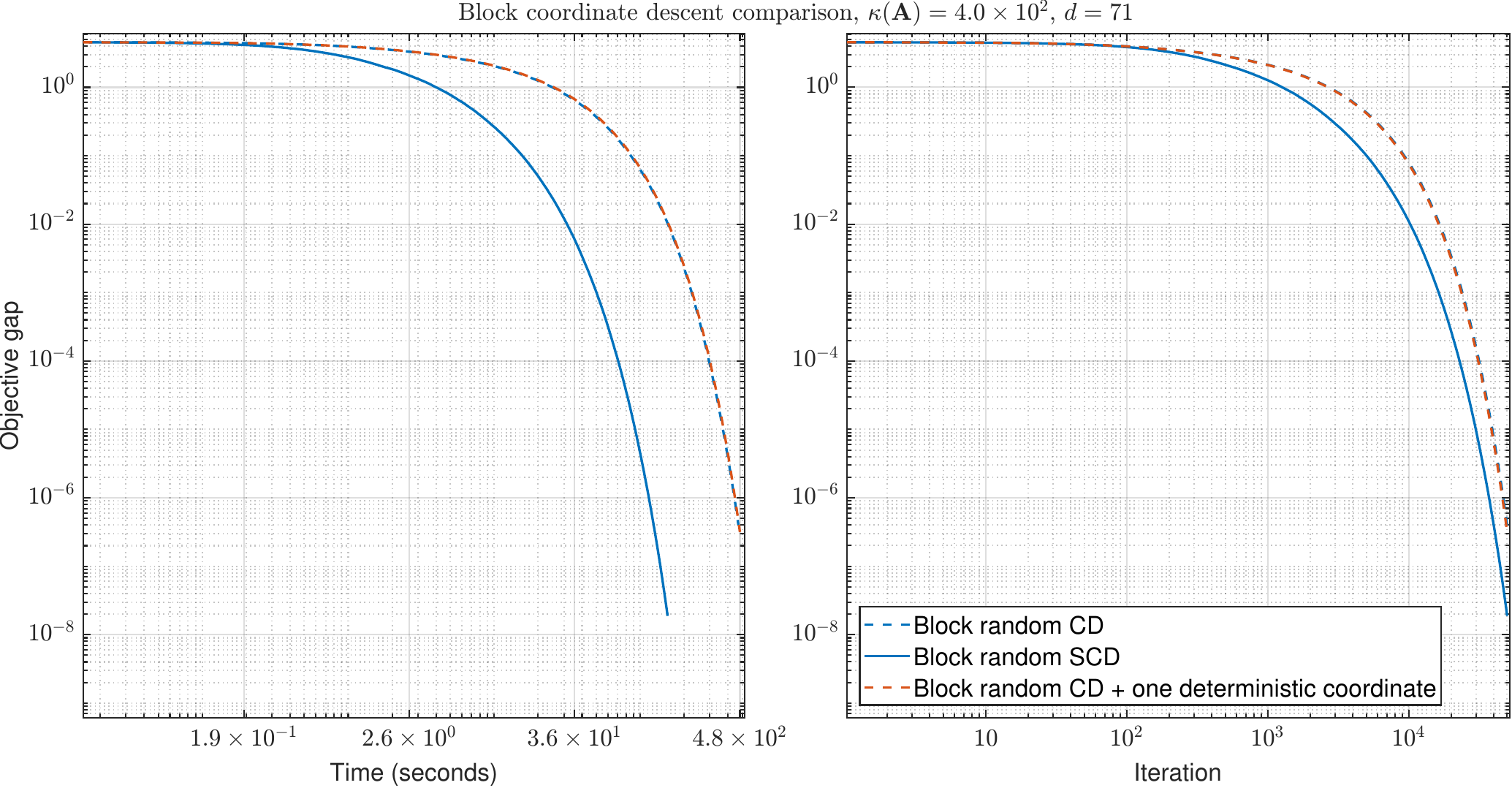}

    %\vspace{0.4cm}

    \includegraphics[width=0.9\linewidth]{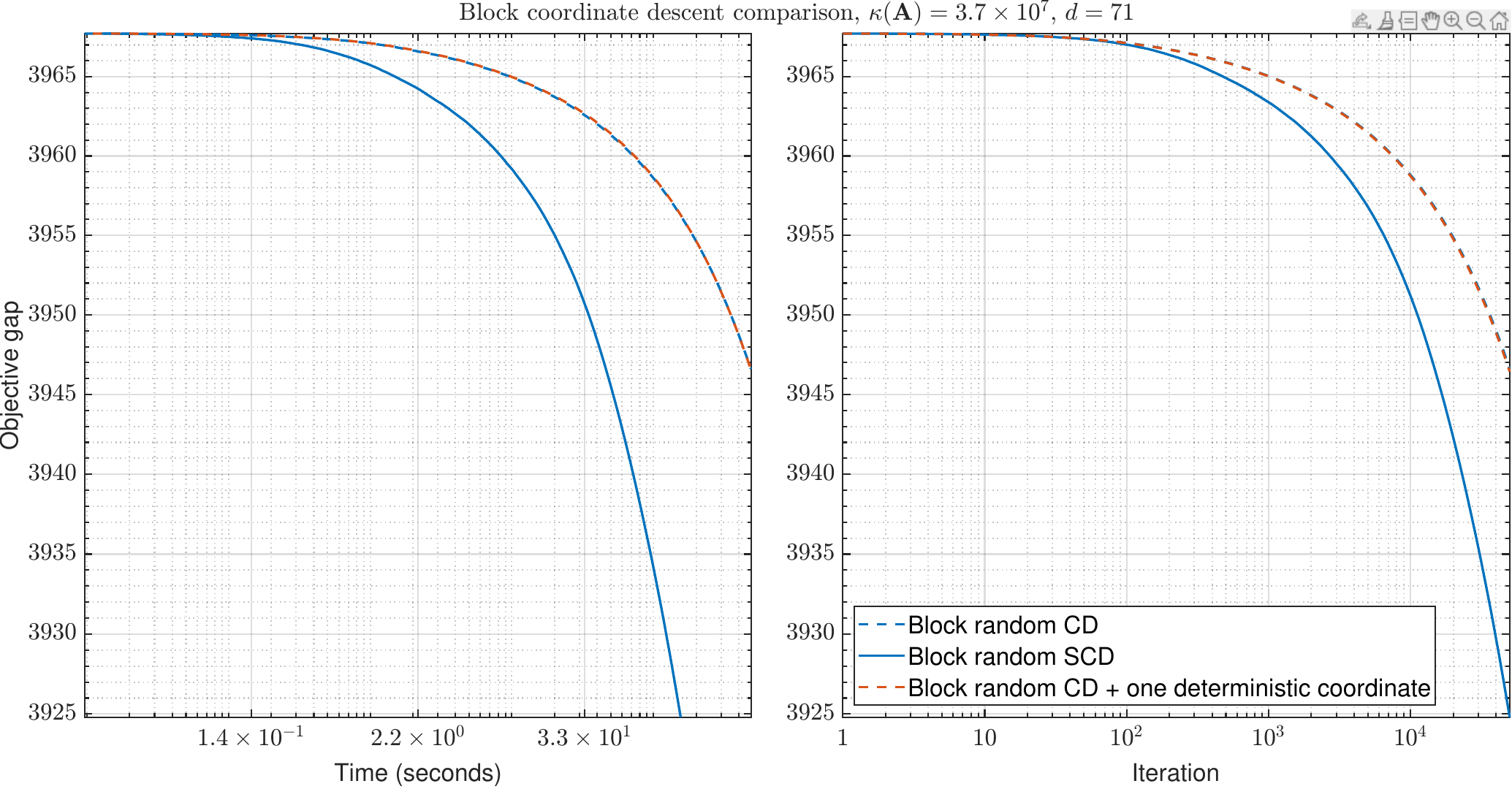}

    \caption{Comparison of the block coordinate descent methods for two test problems with different condition numbers. Results are averaged over 100 realizations.}\label{fig:block_n5000_comparison}
\end{figure}
\subsection{Block coordinate methods}

Figure~\ref{fig:block_n5000_comparison} compares ordinary and scaled block coordinate descent with block size $d= \lceil \sqrt{n} \rceil = 71$ and an iteration budget of $K=50\,000$. The average reported condition numbers are $\kappa(\mathbf{A})\approx4.022\times10^{2}$ in the upper plot and $\kappa(\mathbf{A})\approx3.6977\times10^{7}$ in the lower plot. Results are again averaged over $100$ realizations. The smaller condition number block configuration differs from its single-coordinate counterpart; the two figures should therefore not be interpreted as a controlled comparison of block sizes on identical problems.

At a common block size, BSCD augments the selected coordinate subspace with the current iterate rather than selecting an additional coordinate. In the nondegenerate exact branch, this gives a $72$-dimensional search space while retaining the same $71\times71$ principal-block factorization. The additional right-hand side, coefficient calculations, and global scaling account for its extra work. This comparison examines whether the scaling direction remains useful when the ordinary method already optimizes over a substantial coordinate block. The tests compare the proposed BSCD against the $d$- and $(d+1)$-dimensional ordinary BCD.

Overall, the reported experiments indicate that scaling can improve computational performance, with the benefit particularly apparent in the more ill-conditioned configurations. This is an empirical observation for the tested problems, not a consequence of the condition number alone: the additional decrease depends on the current iterate and its interaction with the selected coordinate or block. The results therefore assess scaling as an inexpensive enlargement of the search space, rather than as a uniformly superior replacement for every coordinate-selection strategy. The effect of the scaling factor is especially striking in this setting, where the $d$- and $ (d+1) $-dimensional ordinary BCD approaches offer nearly identical empirical behavior, but scaling provides significantly better results.

\section{Conclusions}

The proposed method augments an exact coordinate correction with an exact scaling of the current iterate. For symmetric positive definite quadratic objectives, both parameters are obtained from a $2\times2$ problem, and the resulting decrease is governed by a state-dependent effective curvature. The exact scaled update dominates the ordinary one coordinate update for the same selected coordinate, while retaining a linear cost per iteration for dense matrices. The same construction yields adaptive randomized sampling, empirical comparisons with two coordinate methods, and a direct block extension. Empirical experiments show improved results, especially when the working matrices have higher condition numbers.

\section*{Acknowledgment}

This research is supported by the project “Romanian Hub for Artificial Intelligence - HRIA”, Smart Growth, Digitization and Financial Instruments Program, 2021-2027, MySMIS no. 351416.

\bibliography{refs, sample, scaled_coordinate_descent_anm}

\end{document}